\documentclass[11pt]{article}

\usepackage{latexsym}
\usepackage{amssymb}

\newtheorem{theorem}{Theorem}
\newtheorem{proposition}{Proposition}
\newtheorem{corollary}{Corollary}
\newtheorem{definition}{Definition}
\newtheorem{lemma}{Lemma}
\newtheorem{remark}{Remark}

\begin{document}
{

\begin{center}
{\large\bf 
An analytical parameterization for all solutions of the two-dimensional moment problem under Carleman-type conditions}
\end{center}
\begin{center}

{\bf S.M. Zagorodnyuk}
\end{center}

\noindent
\textbf{Abstract.} 
The two-dimensional moment problem consists of finding a positive Borel measure $\mu$ in $\mathbb{R}^2$ such that
$\int_{\mathbb{R}^2} t_1^m t_2^n d\mu = s_{m,n}$, $m,n=0,1,2,...$, where $s_{m,n}$ are prescribed real constants (moments).
We study this moment problem in the case when the sequence $\{ s_{m,n} \}_{m,n=0}^\infty$ is positive semi-definite, and the following Carleman-type
conditions hold: 
$$ \sum_{k=1}^\infty \frac{1}{ \sqrt[2k]{ s_{2m,2k} + s_{2m+2,2k} } } = \infty,\quad  m=0,1,2,.... $$
In this case all solutions of the moment problem are parameterized by a class of analytic contractive operator-valued functions.
The special case of the determinate moment problem is characterized.
We introduce a notion of a generalized resolvent for a pair of commuting symmetric operators.
We use basic properties of such generalized resolvents as a main tool in studying the above moment problem.



\section{Introduction}

As usual, we denote by $\mathbb{R}, \mathbb{C}, \mathbb{N}, \mathbb{Z}, \mathbb{Z}_+$
the sets of real numbers, complex numbers, positive integers, integers and non-negative integers,
respectively. Let $d\in\mathbb{N}$. 
Set $\mathbb{Z}^d_+ = \mathbb{Z}_+\times \ldots \times\mathbb{Z}_+$, and $\mathbb{R}^d = \mathbb{R}\times \ldots \times\mathbb{R}$,
where the Cartesian products are taken with $d$ copies.
For $\mathbf{k} = (k_1,\ldots,k_d)\in\mathbb{Z}^d_+$, $\mathbf{t} = (t_1,\ldots,t_d)\in\mathbb{R}^d$ we denote by
$\mathbf{t}^{\mathbf{k}}$ the monomial $t_1^{k_1}\ldots t_d^{k_d}$.
By $\mathfrak{B}(\mathbb{R}^d)$ we mean the set of all Borel subsets of $\mathbb{R}^d$.

Let $\mathcal{S} = (s_{\mathbf{k}})_{\mathbf{k}\in\mathbb{Z}^d_+}$
be an arbitrary set of real constants (which are called moments).
The (full) multidimensional moment problem consists of finding a (non-negative) measure $\mu$ on $\mathfrak{B}(\mathbb{R}^d)$
such that
\begin{equation}
\label{f1_1}
\int_{ \mathbb{R}^d } \mathbf{t}^{\mathbf{k}} d\mu(\mathbf{t}) = s_{\mathbf{k}},\qquad  \forall \mathbf{k}\in\mathbb{Z}^d_+.
\end{equation}

Various aspects of the multidimensional moment problems were described in the books of Berezansky~\cite{cit_1000_Berezansky_1965__Book},
Berg, Christensen and Ressel~\cite{cit_980_Berg_Christiansen_Ressel__Book}, 
Curto and Fialkow~\cite{cit_985_Curto_Fialkow__Book1}, \cite{cit_985_Curto_Fialkow__Book2},
Lasserre \cite{cit_985_Lasserre_Book},
Schm\"udgen~\cite{cit_1000_Schmudgen_Book_2017}, 
Shohat and Tamarkin~\cite{cit_1000_Shohat_Tamarkin_1943__Book}.
In~\cite{cit_980_Fuglede} Fuglede presented an operator-theoretic approach to the multidimensional moment problem.
The moment problem was transformed into a problem of finding commuting self-adjoint extensions of commuting
densely defined symmetric operators, associated to the given data.
The operator-theoretic approach to moment problems has its origin in papers of Naimark
in 1940--1943 and then it was developed by many authors, see historical notes in~\cite{cit_15000_Zagorodnyuk_2017_J_Adv_Math_Stud}.

For all bounded and some unbounded semi-algebraic sets the existence of solutions of the corresponding $K$-moment problems 
was established, see~\cite{cit_1000_Schmudgen_Book_2017}.
In the case of no prescribed restrictions on the measure support various useful conditions of solvability for the moment problem~(\ref{f1_1})
are known, see discussions on this topic in Berezansky's book~\cite{cit_980_Berg_Christiansen_Ressel__Book}
and in Schm\"udgen's book~\cite{cit_1000_Schmudgen_Book_2017}.
These conditions include Carleman-type conditions which we shall need in the present work. 
Some complicated conditions for the solvability of multidimentional moment problems were given 
in~\cite{cit_1000_Shohat_Tamarkin_1943__Book}, \cite{cit_1100_Stochel_Szafraniec__JFA_1998},
\cite{cit_993_P_V__1999}, \cite{cit_14500_Zagorodnyuk_2010_AFA}, 
\cite{cit_14220_Vasilescu}.

The case $d=1$ is well understood and main questions are answered here.
In this paper we shall focus our attention on the two-dimensional moment problem, $d=2$.
For a notational convenience we shall restate the moment problem in this case. 
\textit{The two-dimensional moment problem} consists of finding a (non-negative) measure $\mu$ on $\mathfrak{B}(\mathbb{R}^2)$
such that
\begin{equation}
\label{f1_1_d_2}
\int_{ \mathbb{R}^2 } t_1^m t_2^n d\mu(t_1,t_2) = s_{m,n},\qquad   m,n=0,1,2,...,
\end{equation}
where $\{ s_{m,n} \}_{m,n=0}^\infty$ are prescribed real constants (called moments).
The following Carleman-type condition will be our starting point in studying this moment problem:
\begin{equation}
\label{f1_2}
\sum_{k=1}^\infty \frac{1}{ \sqrt[2k]{ s_{2m,2k} + s_{2m+2,2k} } } = \infty,\quad \forall m\in\mathbb{Z}_+.
\end{equation}
Notice that this condition is close to Corollary~14.13 in~\cite[p. 365]{cit_1000_Schmudgen_Book_2017}, In fact, one can 
apply the usual Carleman condition to the one-dimensional sequence of Corollary~14.13, and then
interchange the roles of the first and second coordinates.
Our main aim here is to get an analytical parameterization of all solutions to the moment problem~(\ref{f1_1_d_2})
under conditions~(\ref{f1_2}).

Besides the one-dimensional case, in some cases of multidinensional moment problems an analytical parameterization of
all solutions was constructed, see~\cite{cit_14558_Zagorodnyuk_2012_JOT} for Devinatz's moment problem,
and~\cite{cit_14700_Zagorodnyuk_2013_MFAT} for the two-dimensional moment problem in a strip.
In the above papers the operator approach was used. However, generalized resolvents for pairs of operators
and their properties were not available at that time. Instead of that way, some algebraic manipulations
with integral sums were used to obtain an analytical parameterizations of all solutions for 
the above-mentioned moment problems.

Generalized resolvents of symmetric and isometric operators were often used to solve various interpolation problems.
We shall need basic properties of generalized resolvents for single symmetric and isometric operators, see
a survey in~\cite{cit_14550_Zagorodnyuk_2013_Survey__Generalized resolvents}.
A notion of a generalized resolvent for two commuting isometric operators was given
in~\cite{cit_15000_Zagorodnyuk_2017_NYJM_Generalized_resolvents}.
In the present paper we shall introduce a notion of a generalized resolvent for two commuting (not necessarily densely defined)
closed symmetric operators. In the case, when one operator is self-adjoint, an analytical parameterization of
all generalized resolvents is obtained. This parameterization is used to get a parameterization of
all solutions to the moment problem~(\ref{f1_1_d_2}), with a positive semi-definite set of moments and satisfying condition~(\ref{f1_2}).

We notice that Mor\'an in~\cite{cit_991__Moran_1993} described
minimal commuting self-adjoint extensions of two commuting self-adjoint and symmetric operators. Probably, by a symmetric
operator it was meant a densely defined symmetric operator. Otherwise, the argument on page~168, line~4,
of~\cite{cit_991__Moran_1993} is not clear.

The paper is organized as follows. In Section~\ref{s_1} we shall recall some known facts on generalized resolvents and 
obtain some auxiliary results. Then we
introduce a definition for a generalized resolvent of two commuting
closed symmetric operators in a Hilbert space $H$.
In the case when one operator is self-adjoint, we give an
analytical parameterization of all generalized resolvents (Theorem~\ref{t4_2}).
In Section~\ref{s_3} we recall the operator theoretic approach for the two-dimensional moment problem
and relate it to generalized resolvents.
We provide a simple description of canonical solutions (Theorem~\ref{t6_2}).
Then we shall obtain a parameterization of all solutions of the moment problem~(\ref{f1_1_d_2}), when the moment
sequence is positive semi-definite and condition~(\ref{f1_2}) holds (Theorem~\ref{t6_3}). 
The special case of the determinate moment problem is characterized as a corollary (Corollary~\ref{c6_1}).

\noindent
{\bf Notations. }
Besides the given above notations we shall use the following conventions.
We denote $\mathbb{D} = \{ z\in\mathbb{C}\ \arrowvert \ |z| < 1 \}$,  $\mathbb{T} = \{ z\in\mathbb{C}\ \arrowvert \ |z|=1 \}$,
$\mathbb{R}_e = \mathbb{C}\backslash\mathbb{R}$, $\mathbb{T}_e = \mathbb{C}\backslash\mathbb{T}$,
$\mathbb{C}_+ =  \{ z\in \mathbb{C}\ \arrowvert \ \mathop{\rm Im}\nolimits z > 0 \}$.

All Hilbert spaces appearing in this paper are assumed to be separable.
If H is a Hilbert space then $(\cdot,\cdot)_H$ and $\| \cdot \|_H$ mean
the scalar product and the norm in $H$, respectively.
Indices may be omitted in obvious cases.
For a linear operator $A$ in $H$, we denote by $D(A)$
its  domain, by $R(A)$ its range, and $A^*$ means the adjoint operator
if it exists. If $A$ is invertible then $A^{-1}$ means its
inverse. $\overline{A}$ means the closure of the operator, if the
operator is closable. If $A$ is bounded then $\| A \|$ denotes its
norm.
For a set $M\subseteq H$
we denote by $\overline{M}$ the closure of $M$ in the norm of $H$.
By $\mathop{\rm Lin}\nolimits M$ we mean
the set of all linear combinations of elements from $M$,
and $\mathop{\rm span}\nolimits M:= \overline{ \mathop{\rm Lin}\nolimits M }$.
By $A|_M$ we mean the restriction of $A$ to the set $M$.
By $E_H$ we denote the identity operator in $H$, i.e. $E_H x = x$,
$x\in H$. In obvious cases we may omit the index $H$. If $H_1$ is a subspace of $H$, then $P_{H_1} =
P_{H_1}^{H}$ is an operator of the orthogonal projection on $H_1$
in $H$.

If $V$ is a closed isometric operator in $H$, then
$$ M_\zeta = M_\zeta(V) = (E_H - \zeta V) D(V),\quad
N_\zeta = N_\zeta(V) = H\ominus M_\zeta(V), $$ 
where $\zeta\in \mathbb{C}$. We denote
$M_\infty = M_\infty(V) = R(V)$,
$N_\infty = N_\infty(V) = H\ominus R(V)$, and
$$ \mathcal{R}_\zeta = \mathcal{R}_\zeta(V) = (E_H-\zeta V)^{-1},\quad \zeta\in \mathbb{C}\backslash \mathbb{T}. $$

If $A$ is a closed (not necessarily densely defined) symmetric operator in $H$, then
$$ \mathcal{M}_z = \mathcal{M}_z(A) = (A - zE_H) D(A),\quad 
\mathcal{N}_z = \mathcal{N}_z(A) = H\ominus \mathcal{M}_z, $$ 
where $z\in \mathbb{C}$. Moreover, we set
$$ R_z = R_z(A) = (A-zE_H)^{-1},\qquad z\in \mathbb{C}\backslash \mathbb{R}. $$

For a unitary operator $U$ in a Hilbert space $H$ we shall use the following notation:
$$ U(z) := (E_{H} + z U) (E_{H} - z U)^{-1} = -E_{H} + 2 \mathcal{R}_{z}(U),\qquad z\in\mathbb{T}_e. $$

By $\mathcal{S}(D;N,N')$ we denote a class of all analytic in a domain $D\subseteq \mathbb{C}$
operator-valued functions $F(z)$, which values are linear non-expanding operators mapping the whole
$N$ into $N'$, where $N$ and $N'$ are some Hilbert spaces.

\section{Generalized resolvents of commuting symmetric operators}
\label{s_1}

\subsection{Generalized resolvents of single isometric and symmetric operators}

In this subsection we first recall some definitions and facts on generalized resolvents of isometric and symmetric operators 
which we shall need in the sequel.
Then we characterize those generalized resolvents of a closed isometric operator with no non-zero fixed elements, 
which are generated by unitary extensions without nonzero fixed elements (Proposition~\ref{p2_1}).

Let $V$ be a closed isometric operator in a Hilbert space $H$. Consider an arbitrary unitary operator $U\supseteq V$ in
a Hilbert space $\widetilde H\supseteq H$.
The following operator-valued function
\begin{equation}
\label{f2_5}
\mathbf{R}_{u;z} = \mathbf{R}_{u;z}(V) 
= \left.P^{\widetilde H}_H \left(
E_{\widetilde H} - z U 
\right)^{-1}\right|_H,\qquad z\in \mathbb{T}_e,
\end{equation}
is said to be \textit{the generalized resolvent of $V$ (corresponding to the unitary extension $U$)}.
Of course, two different unitary extensions of $V$ can produce the same generalized resolvent $\mathbf{R}_{u;z}(V)$.
It is often useful to consider {minimal unitary extensions. 
Set
\begin{equation}
\label{f2_7}
H_{\mathrm{min}} := \bigvee\limits_{k=-\infty}^\infty U^k H;\quad  U_{\mathrm{min}} := U|_{ H_{\mathrm{min}} },
\end{equation}
where $\bigvee\limits_{k=-\infty}^\infty M_j$ means the closure of a set of all linear combinations of vectors from sets
$M_j\subseteq H$.
The operator $U_{\mathrm{min}}$ is said to be \textit{a minimal unitary extension of $V$ (corresponding to $U$)}.
Properties of minimal unitary extensions are well-known, but it is not easy to
provide a precise reference. For convenience of reading, we shall collect some their relevant properties in the following lemma.

\begin{lemma}
\label{l1_1}
Let $V$ and $U$ be defined as above. The following statements hold:
\begin{itemize}
\item[(i)] The subspace $H_{\mathrm{min}}$ reduces $U$. The operator $U_{\mathrm{min}}$ is a unitary operator in $H_{\mathrm{min}}$;

\item[(ii)] If $W$ is a unitary operator such that $V\subseteq W\subseteq U$, then $W\supseteq U_{\mathrm{min}}$;

\item[(iii)] Suppose that there are two minimal unitary extensions $G_1$, $G_2$ of $V$, acting in some Hilbert
spaces $\mathcal{H}_1$, $\mathcal{H}_2$, respectively. If $G_1$ and $G_2$ produce by~(\ref{f2_5}) the same
generalized resolvent $\mathbf{R}_{u;z}(V)$, then $G_1$ and $G_2$ are unitarily equivalent, i.e., there exists
a unitary operator $K$ which maps $\mathcal{H}_1$ onto the whole $\mathcal{H}_2$, and $G_2 = K G_1 K^{-1}$.
In this case it also holds: $K|_H = E_H$.

\end{itemize}
\end{lemma}

\noindent
\textbf{Proof.} Properties $(i)$-$(ii)$ are obvious. Let us check the validity of property~$(iii)$. Denote
$\mathcal{D}_k = \mathop{\rm Lin}\nolimits \left\{ G_k^j H,\ j\in\mathbb{Z} \right\}$, $k=1,2$. 
Consider the following operator $K_0$ which maps $\mathcal{D}_1$ onto $\mathcal{D}_2$:
\begin{equation}
\label{f2_9}
K_0 \sum_{j=-\infty}^\infty \alpha_j G_1^j h_j = \sum_{j=-\infty}^\infty \alpha_j G_2^j h_j,\quad \alpha_j\in\mathbb{C},\ h_j\in H,
\end{equation}
where all but a finite number of $\alpha_j$s are zeros.
Let us check that $K_0$ is well-defined. Suppose that $g\in\mathcal{D}_1$ has two representations:
\begin{equation}
\label{f2_12}
g = \sum_{j=-\infty}^\infty \alpha_j G_1^j h_j = \sum_{j=-\infty}^\infty \beta_j G_1^j u_j,\qquad \alpha_j,\beta_j\in\mathbb{C},\ h_j,u_j\in H.
\end{equation}
Then
$$ \left\|
\sum_{j=-\infty}^\infty \alpha_j G_2^j h_j - \sum_{j=-\infty}^\infty \beta_j G_2^j u_j
\right\|_{\mathcal{H}_2}^2 = $$
$$ = \left(
\sum_{j=-\infty}^\infty G_2^j ( \alpha_j h_j - \beta_j u_j ),
\sum_{k=-\infty}^\infty G_2^k ( \alpha_k h_k - \beta_k u_k )
\right)_{\mathcal{H}_2} = $$
$$ = 
\sum_{j,k = -\infty}^\infty
\left(
G_2^{j-k} ( \alpha_j h_j - \beta_j u_j ),
\alpha_k h_k - \beta_k u_k 
\right)_{\mathcal{H}_2} = $$
$$ = 
\sum_{j,k = -\infty}^\infty
\int_\mathbb{T}
z^{j-k} d \left(
\mathbf{F}_2
( \alpha_j h_j - \beta_j u_j ), \alpha_k h_k - \beta_k u_k 
\right)_{\mathcal{H}_2} = $$
$$ = 
\sum_{j,k = -\infty}^\infty
\int_\mathbb{T}
z^{j-k} d \left(
\mathbf{F}_1
( \alpha_j h_j - \beta_j u_j ), \alpha_k h_k - \beta_k u_k 
\right) = $$
$$ = \left\|
\sum_{j=-\infty}^\infty \alpha_j G_1^j h_j - \sum_{j=-\infty}^\infty \beta_j G_1^j u_j
\right\|^2 = 0, $$
where $\mathbf{F}_k$ denote the spectral measure of $V$ corresponding to $G_k$.
Here we have used the assumption that the corresponding generalized resolvents coincide (and therefore the spectral functions are the same).
Thus, $K_0$ is well-defined. For arbitrary $h_j, g_j\in H$ we may write
$$ \left(
K_0 \sum_{j=-\infty}^\infty G_1^j h_j, K_0 \sum_{k=-\infty}^\infty G_1^k g_k
\right) = \left(
\sum_{j=-\infty}^\infty G_2^j h_j, 
\sum_{k=-\infty}^\infty G_2^k g_k \right) = $$
$$ = \sum_{j,k = -\infty}^\infty
\left(
G_2^{j-k} h_j, g_k \right) = 
\int_\mathbb{T} z^{j-k} d(\mathbf{F}_2 h_j, g_k) = \int_\mathbb{T} z^{j-k} d(\mathbf{F}_1 h_j, g_k) =  
$$
$$ = \left(
\sum_{j=-\infty}^\infty G_1^j h_j, \sum_{k=-\infty}^\infty G_1^k g_k
\right). $$
Thus, $K_0$ is isometric. By the continuity we extend $K_0$ to a unitary operator $K$ which maps $\mathcal{H}_1$ onto $\mathcal{H}_2$.
It is readily checked that $G_2 = K G_1 K^{-1}$ and $K|_H = E_H$. $\Box$

By Chumakin's formula (see, e.g.,~Theorem~2.7 in~\cite[p. 191]{cit_14550_Zagorodnyuk_2013_Survey__Generalized resolvents}), 
an arbitrary generalized resolvent $\mathbf{R}_{u;z}$ has the following form:
\begin{equation}
\label{f2_15}
\mathbf{R}_{u;z} = \left[ E_H - z (V\oplus \Phi_z) \right]^{-1},\qquad z\in \mathbb{D},
\end{equation}
where $\Phi_z\in\mathcal{S}(\mathbb{D};N_0(V),N_\infty(V))$ (see Notations).
Conversely, an arbitrary function $\Phi_z\in\mathcal{S}(\mathbb{D};N_0(V),N_\infty(V))$
defines by relation~(\ref{f2_15}) a generalized resolvent $\mathbf{R}_{u;z}$ of $V$.
Moreover, this correspondence is one-to-one.


Let $A$ be an arbitrary closed symmetric operator in a Hilbert space $H$.
The domain of $A$ is not supposed to be necessarily dense in $H$.
Consider an arbitrary self-adjoint extension $\widehat A$ of $A$, acting in a Hilbert space  $\widehat H\supseteq H$.
The following operator-valued function:
$$ \mathbf{R}_{s;\lambda} = \mathbf{R}_{s;\lambda}(A) =
\left. P^{\widehat{H}}_H \left( \widehat A - \lambda E_{\widehat{H}} \right)^{-1} \right|_H,\qquad
\lambda\in \mathbb{R}_e, $$
is said to be \textit{a generalized resolvent of $A$ (corresponding to $\widehat A$)}.
Consider the Cayley transformation of the operator $A$:
\begin{equation}
\label{f2_17}
U_i = U_i(A) = (A + i E_H)(A - i E_H)^{-1} = E_H + 2i (A - i E_H)^{-1}.
\end{equation}
Then the following operator:
\begin{equation}
\label{f2_19}
W_i := (\widehat A + i E_{\widehat H})(\widehat A - i E_{\widehat H})^{-1}
= E_{\widehat H} + 2i (\widehat A - i E_{\widehat H})^{-1}
\end{equation}
is a unitary extension of the operator $U_i$, having no non-zero fixed points.
The following relation (see~Theorem~2.9 in~\cite[p. 199]{cit_14550_Zagorodnyuk_2013_Survey__Generalized resolvents}):
\begin{equation}
\label{f2_25}
\mathbf{R}_{u;\frac{\lambda - i}{\lambda + i }} (U_i)
= \frac{\lambda + i}{ 2i} E_H +
\frac{(\lambda + i)(\lambda - i)}{ 2i} \mathbf{R}_{s;\lambda} (A),\qquad
\lambda\in \mathbb{R}_e\backslash\{ \pm i \},
\end{equation}
establishes a one-to-one correspondence between all generalized resolvents $\mathbf{R}_{s;\lambda}(A)$
of $A$ and those generalized resolvents $\mathbf{R}_{u;z} (U_i)$ of the closed isometric
operator $U_i$ which are generated by unitary extensions of $U_i$ without non-zero fixed elements.
In the case $\overline{ D(A) } = H$, equality~(\ref{f2_25})
establishes a one-to-one correspondence between all generalized resolvents $\mathbf{R}_{s;\lambda}(A)$
of $A$ and all generalized resolvents $\mathbf{R}_{u;z} (U_i)$ of $U_i$.

In order to formulate the Shtraus formula for the generalized resolvents of $A$, we need to recall some definitions.
Consider the following operator:
$$ X_i \psi = X_i(A) \psi
= \varphi,\qquad \psi\in \mathcal{N}_i(A)\cap \left( \mathcal{N}_{-i}(A) \dotplus D(A)
\right), $$
where $\varphi\in \mathcal{N}_{-i}(A)$: $\psi - \varphi\in D(A)$.
The operator $X_i=X_i(A)$ is said to be \textit{forbidden with respect to the symmetric operator $A$}.
Denote
$\mathbb{C}_+^\varepsilon = \{ z\in\mathbb{C}_+\ \arrowvert\  \varepsilon < \arg z < \pi - \varepsilon \}$,
$0 < \varepsilon < \frac{\pi}{2}$.
A function $F(\lambda)\in \mathcal{S}(\mathbb{C}_+; \mathcal{N}_{i}(A), \mathcal{N}_{ -i }(A))$ 
is said to be \textit{admissible with respect to the operator $A$},
if the validity of
\begin{equation}
\label{AFA_2013___f11_1_p2_1}
\lim_{\lambda\in\mathbb{C}_+^\varepsilon,\ \lambda\to\infty} F(\lambda) \psi = X_{i}\psi,
\end{equation}
\begin{equation}
\label{AFA_2013___f11_2_p2_1}
\underline{\lim}_{\lambda\in\mathbb{C}_+^\varepsilon,\ \lambda\to\infty}
\left[
|\lambda| (\| \psi \|_H - \| F(\lambda) \psi \|_H)
\right] < +\infty,
\end{equation}
for some $\varepsilon$: $0<\varepsilon <\frac{\pi}{2}$,
implies $\psi = 0$.

A set of all operator-valued functions $F(\lambda)\in \mathcal{S}(\mathbb{C}_+; \mathcal{N}_{i}(A),
\mathcal{N}_{ -i }(A))$, which are admissible with respect to the operator $A$,
we shall denote by $\mathcal{S}_{a} (\mathbb{C}_+; \mathcal{N}_{i}(A), \mathcal{N}_{ -i }(A))$.
In the case $\overline{D(A)} = H$, we have 
$\mathcal{S}_{a} (\mathbb{C}_+; \mathcal{N}_{i}(A), \mathcal{N}_{ -i }(A)) = 
\mathcal{S} (\mathbb{C}_+; \mathcal{N}_{i}(A), \mathcal{N}_{ -i }(A))$.

Let $F$ be an arbitrary linear bounded operator with the domain
$D(F) = \mathcal{N}_i(A)$ and range $R(F) \subseteq \mathcal{N}_{-i} (A)$.
Set $U_{i;F} = U_{i;F}(A) = U_i \oplus F$. 
Suppose that $U_{i;F}(A)$ has no non-zero fixed elements. 
In this case we may define the following linear operator in $H$:
$$ A_F = A_{F,i} =
(i U_{i;F} + i E_H) (U_{i;F} - E_H)^{-1}. $$
Notice that the operator $A_{F}$ is an extension of the operator $A$.
The operator $A_F=A_{F,i}$ is said to be \textit{a quasi-self-adjoint extension of $A$,
defined by the operator $F$}.

An arbitrary generalized resolvent $\mathbf{R}_{s;\lambda}$ of $A$ has the following form:
\begin{equation}
\label{AFA_2013___f11_3_p2_1}
\mathbf R_{s;\lambda} = \left\{ \begin{array}{cc}
\left( A_{F(\lambda)} - \lambda E_H \right)^{-1}, & \lambda\in \mathbb{C}_+ \\
\left( A_{F^*(\overline{\lambda})} - \lambda E_H \right)^{-1}, &
\overline{\lambda}\in \mathbb{C}_+
\end{array}
\right.,
\end{equation}
where $F(\lambda)$ is a function from $\mathcal{S}_{a} (\mathbb{C}_+; \mathcal{N}_{i}(A), \mathcal{N}_{ -i }(A))$.
Conversely, an arbitrary function $F(\lambda)\in \mathcal{S}_{a} (\mathbb{C}_+; \mathcal{N}_{i}(A), \mathcal{N}_{ -i }(A))$
defines by relation~(\ref{AFA_2013___f11_3_p2_1}) a generalized resolvent
$\mathbf{R}_{s;\lambda}$ of $A$. This correspondence is a bijection.
Formula~(\ref{AFA_2013___f11_3_p2_1}) was obtained by Shtraus and it is called the Shtraus formula for generalized resolvents
of a symmetric operator. For proofs and more details on this theory we refer to 
a survey~\cite{cit_14550_Zagorodnyuk_2013_Survey__Generalized resolvents}.

Let us now return to the correspondence~(\ref{f2_25}). By the Shtraus formula, for each $F(\lambda)$ from 
$\mathcal{S}_{a} (\mathbb{C}_+; \mathcal{N}_{i}(A), \mathcal{N}_{ -i }(A))$ there corresponds a generalized
resolvent $\mathbf{R}_{s;\lambda}(A)$. By the correspondence~(\ref{f2_25}) for $\mathbf{R}_{s;\lambda}(A)$ there corresponds
a generalized resolvent $\mathbf{R}_{u;z} (U_i)$ of the operator $U_i$. Finally, by Chumakin's formula for 
$\mathbf{R}_{u;z} (U_i)$ there corresponds a function 
$\Phi_z\in\mathcal{S}(\mathbb{D};N_0(U_i),N_\infty(U_i)) = \mathcal{S} (\mathbb{D}; \mathcal{N}_{i}(A), \mathcal{N}_{ -i }(A))$.
Functions $F(\lambda)$ and $\Phi_z$ are related in the following way (see~formula~(4.87) 
in~\cite[p. 281]{cit_14550_Zagorodnyuk_2013_Survey__Generalized resolvents}):
\begin{equation}
\label{AFA_2013___f7_7_p3_1}
F(\lambda) = \Phi_{ \frac{\lambda - i}{\lambda + i} },\qquad \lambda\in \mathbb{C}_+.
\end{equation}
Therefore $\Phi_{ \frac{\lambda - i}{\lambda + i} }\in \mathcal{S}_{a} (\mathbb{C}_+; \mathcal{N}_{i}(A), \mathcal{N}_{ -i }(A))$.

\begin{definition}
\label{d2_1}
Let $A$ be a closed symmetric operator in a Hilbert space $H$, and
$V=U_i$ be its Cayley transformation as in~(\ref{f2_17}).
Denote by $\mathcal{S}_{ \mathrm{adm}; V }(\mathbb{D};N_0(V),N_\infty(V))$ a set of all functions 
$\Phi_z$ from $\mathcal{S}(\mathbb{D};N_0(V),N_\infty(V))$ such that
$$ \Phi_{ \frac{\lambda - i}{\lambda + i} }\in \mathcal{S}_{a} (\mathbb{C}_+; \mathcal{N}_{i}(A), \mathcal{N}_{ -i }(A)). $$
\end{definition}

\begin{remark}
\label{r2_1}
If $\overline{D(A)} = H$, then
$\mathcal{S}_{a} (\mathbb{C}_+; \mathcal{N}_{i}(A), \mathcal{N}_{ -i }(A)) = 
\mathcal{S} (\mathbb{C}_+; \mathcal{N}_{i}(A), \mathcal{N}_{ -i }(A))$, and 
we have $\mathcal{S}_{ \mathrm{adm}; V }(\mathbb{D};N_0(V),N_\infty(V)) = \mathcal{S}(\mathbb{D};N_0(V),N_\infty(V))$.
\end{remark}

\begin{proposition}
\label{p2_1}
Let $V$ be a closed isometric operator in a Hilbert space $H$ with no nonzero fixed elements, and
$A = i ( V + E_H ) (V - E_H)^{-1}$.
A generalized resolvent $\mathbf{R}_{u;z}$ of $V$ is
generated by a unitary extension without nonzero fixed elements if and only if
$\Phi_z\in \mathcal{S}_{ \mathrm{adm}; V }(\mathbb{D};N_0(V),N_\infty(V))$, where
$\Phi_z$ is the parameter in Chumakin's formula~(\ref{f2_15}) for $\mathbf{R}_{u;z}$.
\end{proposition}

\noindent
\textbf{Proof.} \textit{Sufficiency.} Suppose that $\widetilde \Phi_z\in \mathcal{S}_{ \mathrm{adm}; V }(\mathbb{D};N_0(V),N_\infty(V))$.
Then $F(\lambda) := \widetilde \Phi_{ \frac{\lambda - i}{\lambda + i} }$, ($\lambda\in \mathbb{C}_+$), belongs
to $\mathcal{S}_{a} (\mathbb{C}_+; \mathcal{N}_{i}(A), \mathcal{N}_{ -i }(A))$.
Repeating the construction for $F(\lambda)$, given before the statement of the proposition, we get
a generalized resolvent $\mathbf{R}_{u;z} (V)$ of the operator $V=U_i(A)$, which is generated by a unitary
extension without nonzero fixed elements.
By formula~(\ref{AFA_2013___f7_7_p3_1}) we see that $\widetilde \Phi_z$ coincides with the parameter $\Phi_z$ in Chumakin's formula
for $\mathbf{R}_{u;z} (V)$.

\noindent
\textit{Necessity.} Suppose that a generalized resolvent $\mathbf{R}_{u;z}$ of $V$ is
generated by a unitary extension without nonzero fixed elements.
By formula~(\ref{f2_25}) $\mathbf{R}_{u;z}$ is related with a generalized resolvent 
$\mathbf{R}_{s;\lambda}(A)$ of $A$.
Let $F(\lambda)$ be the parameter in the Shtraus formula, corresponding to $\mathbf{R}_{s;\lambda}(A)$.
By the construction for $F(\lambda)$, given before the statement of the proposition, we get
$\Phi_{ \frac{\lambda - i}{\lambda + i} }\in \mathcal{S}_{a} (\mathbb{C}_+; \mathcal{N}_{i}(A), \mathcal{N}_{ -i }(A))$,
i.e. $\Phi_z\in \mathcal{S}_{ \mathrm{adm}; V }(\mathbb{D};N_0(V),N_\infty(V))$.
$\Box$

\subsection{Generalized resolvents for pairs of commuting symmetric operators}

We are going to introduce a notion of a generalized resolvent for a pair of commuting symmetric operators.
One could imagine such a function as a compression of a product of usual resolvents for commuting self-adjoint extensions.
However, the definition will be a little bit different.

Let $A_1$ and $A_2$ be two closed (not necessarily densely defined) symmetric operators in a Hilbert space $H$, such that 
their Cayley's transformations
\begin{equation}
\label{f3_5}
V_j := (A_j + i E_H) (A_j - i E_H)^{-1},\qquad j=1,2, 
\end{equation}
commute whenever possible:
\begin{equation}
\label{f3_7}
V_1 V_2 h = V_2 V_1 h,\qquad  h\in D(V_1 V_2) \cap D(V_2 V_1).
\end{equation}
Suppose that there exist commuting self-adjoint operators $B_1, B_2$ in a Hilbert space
$\widetilde H\supseteq H$ such that $B_1\supseteq A_1$, $B_2\supseteq A_2$. 
The following operator-valued function
$$ \mathbf{R}_{s;\lambda_1,\lambda_2} = \mathbf{R}_{s;\lambda_1,\lambda_2} (A_1, A_2) = $$
$$ =
\left. 
P^{\widetilde H}_H (E_{\widetilde H} + \lambda_1 B_1) (B_1 - \lambda_1 E_{\widetilde H})^{-1}
(E_{\widetilde H} + \lambda_2 B_2) (B_2 - \lambda_2 E_{\widetilde H})^{-1}
\right|_H, $$
\begin{equation}
\label{f3_10}
\lambda_1, \lambda_2\in \mathbb{R}_e,
\end{equation}
is said to be \textit{a generalized resolvent of a pair of symmetric operators $A_1, A_2$ (corresponding to extensions $B_1, B_2$)}.
Notice that it satisfies the following relation, which we shall need later:
\begin{equation}
\label{f3_11}
\mathbf{R}_{s;\lambda_1,\lambda_2}^* = \mathbf{R}_{s;\overline{ \lambda_1 }, \overline{ \lambda_2 }},\qquad
\lambda_1, \lambda_2\in \mathbb{R}_e.
\end{equation}

Consider Cayley's transformations
\begin{equation}
\label{f3_12}
U_j := (B_j + i E_{\widetilde H}) (B_j - i E_{\widetilde H})^{-1},\qquad j=1,2, 
\end{equation}
which are commuting unitary extensions of isometric operators $V_1,V_2$.
Recall that the following function (see~\cite{cit_15000_Zagorodnyuk_2017_NYJM_Generalized_resolvents})
$$ \mathbf{R}_{u;z_1,z_2} = \mathbf{R}_{u;z_1,z_2} (V_1,V_2) = $$
$$ =
\left. 
P^{\widetilde H}_H (E_{\widetilde H} + z_1 U_1) (E_{\widetilde H} - z_1 U_1)^{-1}
(E_{\widetilde H} + z_2 U_2) (E_{\widetilde H} - z_2 U_2)^{-1} \right|_H, $$
\begin{equation}
\label{f3_15}
z_1, z_2\in \mathbb{T}_e,
\end{equation}
is said to be \textit{a generalized resolvent of a pair of isometric operators $V_1, V_2$ (corresponding to extensions $U_1,U_2$)}.
We are going to establish a correspondence between $\mathbf{R}_{s;\lambda_1,\lambda_2}(A_1, A_2)$ and 
$\mathbf{R}_{u;z_1,z_2} (V_1,V_2)$, similar to formula~(\ref{f2_25}).

\begin{proposition}
\label{p2_2}
Let $A_1,A_2$ be two closed (not necessarily densely defined) symmetric operators in a Hilbert space $H$, such that 
their Cayley's transformations $V_j := (A_j + i E_H) (A_j - i E_H)^{-1}$, $j=1,2$, satisfy relation~(\ref{f3_7}).
Suppose that there exists a generalized resolvent of symmetric operators $A_1,A_2$.
The following relation:
\begin{equation}
\label{f3_25}
\mathbf{R}_{u; \frac{\lambda_1 - i}{\lambda_1 + i }, \frac{\lambda_2 - i}{\lambda_2 + i }} (V_1,V_2)
= - \mathbf{R}_{s;\lambda_1,\lambda_2} (A_1,A_2),\qquad
\lambda_1,\lambda_2\in \mathbb{R}_e\backslash\{ \pm i \},
\end{equation}
establishes a one-to-one correspondence between all generalized resolvents $\mathbf{R}_{s;\lambda_1,\lambda_2}(A_1,A_2)$
of $A_1,A_2$ and those generalized resolvents $\mathbf{R}_{u;z_1,z_2} (V_1,V_2)$ of closed isometric
operators $V_1,V_2$ which are generated by commuting extensions of $V_1,V_2$ without non-zero fixed elements.
 
\end{proposition}

\noindent
\textbf{Proof.} Let $\mathbf{R}_{s;\lambda_1,\lambda_2}(A_1,A_2)$ be an arbitrary generalized resolvent of $A_1,A_2$.
Then $\mathbf{R}_{s;\lambda_1,\lambda_2}(A_1,A_2)$ has a representation of type~(\ref{f3_10}), generated by
some commuting self-adjoint extensions $B_1, B_2$ in a Hilbert space $\widetilde H\supseteq H$.
Define $U_1,U_2$ as in~(\ref{f3_12}). Of course, they have no nonzero fixed points.
Choose arbitrary $\lambda_1,\lambda_2\in \mathbb{R}_e\backslash\{ \pm i \}$, and
set $z_j = (\lambda_j - i)/(\lambda_j +i)$, $j=1,2$. Then $z_1,z_2\in \mathbb{T}_e\backslash\{ 0 \}$.
We may write
$$ R_{z_1,z_2}(U_1,U_2) :=  (E_{\widetilde H} + z_1 U_1) (E_{\widetilde H} - z_1 U_1)^{-1}
(E_{\widetilde H} + z_2 U_2) (E_{\widetilde H} - z_2 U_2)^{-1} = $$
$$ = 
\left(
- E_{\widetilde H} + 2 (E_{\widetilde H} - z_1 U_1)^{-1}
\right)
\left(
- E_{\widetilde H} + 2 (E_{\widetilde H} - z_2 U_2)^{-1}
\right). $$
Observe that (cf.~\cite[p. 199]{cit_15000_Zagorodnyuk_2017_NYJM_Generalized_resolvents})
$$ (E_{\widetilde H} - z_j U_j)^{-1} = 
\frac{ \lambda_j + i }{ 2i } E_{\widetilde H} +
\frac{ (\lambda_j + i)(\lambda_j - i) }{2i} (B_j - \lambda_j E_{\widetilde H})^{-1},\qquad
j=1,2. $$
Then
$$ R_{z_1,z_2}(U_1,U_2) = -
(E_{\widetilde H} + \lambda_1 B_1) (B_1 - \lambda_1 E_{\widetilde H})^{-1}
(E_{\widetilde H} + \lambda_2 B_2) (B_2 - \lambda_2 E_{\widetilde H})^{-1}. $$
Applying the projection operator and restricting to $H$, we obtain relation~(\ref{f3_25}).
Notice that the generalized resolvent $\mathbf{R}_{u;z_1,z_2} (V_1,V_2)$ can be defined by~(\ref{f3_25})
for all points $z_1,z_2\in \mathbb{T}_e\backslash\{ 0 \}$. At remaining points $\mathbf{R}_{u;z_1,z_2} (V_1,V_2)$ 
can be determined by the weak continuity.

Conversely, choose an arbitrary generalized resolvent $\mathbf{R}_{u;z_1,z_2} (V_1,V_2)$, having a representation of type~(\ref{f3_15})
where $U_1,U_2$ have no nonzero fixed elements. We may define self-adjoint operators $B_1, B_2$ by relation~(\ref{f3_12}). 
Let $\mathbf{R}_{s;\lambda_1,\lambda_2}(A_1,A_2)$ be the generalized resolvent of $A_1,A_2$, generated by extensions $B_1, B_2$.
Repeating the construction at the beginning of the proof we shall obtain formula~(\ref{f3_25}).

If the left-hand side of~(\ref{f3_25}) gives the same function for two generalized resolvents $\mathbf{R}'_{s;\lambda_1,\lambda_2}(A_1,A_2)$
and $\mathbf{R}''_{s;\lambda_1,\lambda_2}(A_1,A_2)$ then
$\mathbf{R}'_{s;\lambda_1,\lambda_2}(A_1,A_2) =\mathbf{R}''_{s;\lambda_1,\lambda_2}(A_1,A_2)$ for
$\lambda_1,\lambda_2\in \mathbb{R}_e\backslash\{ \pm i \}$. For the remaining points the equality follows by the weak continuity.
$\Box$

From the formulation of Proposition~\ref{p2_2} there appear the following natural question:

\noindent
\textbf{Question}: \textit{how one can check that $\mathbf{R}_{u;z_1,z_2} (V_1,V_2)$ is generated by 
commuting extensions of $V_1,V_2$ without non-zero fixed elements?}
An answer will be given in Theorem~\ref{t3_1} below.

Let $\mathbf{R}_{u;z_1,z_2} (V_1,V_2)$ be an arbitrary generalized resolvent of $V_1,V_2$. Then a relation of type~(\ref{f3_15})
holds with some unitary operators $U_1,U_2$.
We may write
$$ \mathbf{R}_{u;z_1,0} = - E_H + 2 \left. P^{\widetilde H}_H (E_{\widetilde H} - z_1 U_1)^{-1} \right|_H = 
- E_H + 2 \mathbf{R}_{u;z_1}(V_1),\ z_1\in\mathbb{T}_e, $$
and
$$ \mathbf{R}_{u;0,z_2} = - E_H + 2 \left. P^{\widetilde H}_H (E_{\widetilde H} - z_2 U_2)^{-1} \right|_H = 
- E_H + 2 \mathbf{R}_{u;z_2}(V_2),\ z_2\in\mathbb{T}_e, $$
where
\begin{equation}
\label{f3_26}
\mathbf{R}_{u;z_j}(V_j) := \left. P^{\widetilde H}_H (E_{\widetilde H} - z_j U_j)^{-1} \right|_H,\qquad j=1,2.
\end{equation}
Notice that
\begin{equation}
\label{f3_26_2}
\mathbf{R}_{u;z_1}(V_1) = (1/2) (\mathbf{R}_{u;z_1,0} + E_H),\ 
\mathbf{R}_{u;z_2}(V_2) = (1/2) (\mathbf{R}_{u;0,z_2} + E_H).
\end{equation}
So, the values of $\mathbf{R}_{u;z_1}(V_1)$ and $\mathbf{R}_{u;z_2}(V_2)$ do not depend on 
a particular choice of $U_1,U_2$, generating $\mathbf{R}_{u;z_1,z_2} (V_1,V_2)$.
By Chumakin's formula~(\ref{f2_15}) we may write
$$ (1/2) (\mathbf{R}_{u;z_1,0} + E_H) =  \left[ E_H - z_1 (V_1\oplus \Phi_1(z_1) \right]^{-1},\ z_1\in\mathbb{D}, $$
and
$$ (1/2) (\mathbf{R}_{u;0,z_2} + E_H) = \left[ E_H - z_2 (V_2\oplus \Phi_2(z_2) \right]^{-1},\ z_2\in\mathbb{D}, $$
where 
$\Phi_j(z_j)\in\mathcal{S}(\mathbb{D};N_0(V_j),N_\infty(V_j))$, $j=1,2$.
Then
\begin{equation}
\label{f3_27}
\Phi_1(z_1) = \left.\left(
\frac{1}{z_1} E_H - \frac{2}{z_1} (\mathbf{R}_{u;z_1,0} + E_H)^{-1}
\right)\right|_{H\ominus D(V_1)},\qquad z_1\in\mathbb{D}\backslash\{ 0 \}, 
\end{equation}
\begin{equation}
\label{f3_29}
\Phi_2(z_2) = \left.\left(
\frac{1}{z_2} E_H - \frac{2}{z_2} (\mathbf{R}_{u;0,z_2} + E_H)^{-1}
\right)\right|_{H\ominus D(V_2)},\qquad z_2\in\mathbb{D}\backslash\{ 0 \}. 
\end{equation}
If $U_1$ and $U_2$ have no nonzero fixed elements, then
by~Proposition~\ref{p2_1} we get
\begin{equation}
\label{f3_35}
\Phi_j(z_j)\in \mathcal{S}_{ \mathrm{adm}; V_j }(\mathbb{D};N_0(V_j),N_\infty(V_j)),\qquad j=1,2.
\end{equation}

\begin{theorem}
\label{t3_1}
In conditions of Proposition~\ref{p2_2}, let $\mathbf{R}_{u;z_1,z_2} (V_1,V_2)$ be an arbitrary generalized resolvent of
$V_1,V_2$. The generalized resolvent $\mathbf{R}_{u;z_1,z_2} (V_1,V_2)$ is generated by~(\ref{f3_15}) with some commuting unitary extensions 
$U_1, U_2$ without nonzero fixed elements if and only if relation~(\ref{f3_35}) holds, where 
$\Phi_1(z_1), \Phi_2(z_2)$ ($z_1,z_2\in\mathbb{D}$) are defined by formulas~(\ref{f3_27}),(\ref{f3_29}),
using the weak continuity at zero.

\end{theorem}

\noindent
\textbf{Proof.} The necessity was checked before the statement of the theorem. Let us check the sufficiency.
Let $\mathbf{R}_{u;z_1,z_2} (V_1,V_2)$ be a generalized resolvent of $V_1,V_2$, and relation~(\ref{f3_35}) holds.
As we have seen, $\Phi_1(z_1), \Phi_2(z_2)$ are the parameters in Chumakin's formula for the generalized resolvents
$\mathbf{R}_{u;z_1}(V_1)$ and $\mathbf{R}_{u;z_2}(V_2)$ as in~(\ref{f3_26_2}).
By Proposition~\ref{p2_1} condition~(\ref{f3_35}) guarantees that the latter generalized resolvents are generated by some unitary extensions 
without nonzero fixed elements $S_1\supseteq V_1$
and $S_2\supseteq V_2$ in some Hilbert spaces $K_1$ and $K_2$, respectively. 
Moreover,
$\mathbf{R}_{u;z_1}(V_1)$ and $\mathbf{R}_{u;z_2}(V_2)$ are generated by the corresponding
minimal unitary extensions $S_{1,\mathrm{min}}$ and $S_{2,\mathrm{min}}$ (cf.~(\ref{f2_7})).
Of course, restrictions $S_{1,\mathrm{min}}$ and $S_{2,\mathrm{min}}$ have no nonzero fixed elements as well.

By the definition of a generalized resolvent, $\mathbf{R}_{u;z_1,z_2} (V_1,V_2)$ admits the following representation:
$$ \mathbf{R}_{u;z_1,z_2} = \mathbf{R}_{u;z_1,z_2} (V_1,V_2) = $$
$$ =
\left. 
P^{K}_H (E_{K} + z_1 W_1) (E_{K} - z_1 W_1)^{-1}
(E_{K} + z_2 W_2) (E_{K} - z_2 W_2)^{-1} \right|_H, $$
\begin{equation}
\label{f3_45}
z_1, z_2\in \mathbb{T}_e,
\end{equation}
where $W_1$ and $W_2$ are commuting unitary operators in a Hilbert space $K\supseteq H$, and 
$W_1 \supseteq V_1$, $W_2 \supseteq V_2$. 
Denote by $W_{j,\mathrm{min}}$ the minimal unitary extension of $V_j$,
generated by $W_j$, acting on a Hilbert space $K_{j,\mathrm{min}}\supseteq H$, j=1,2.

Notice that $\mathbf{R}_{u;z_1}(V_1)$ and $\mathbf{R}_{u;z_2}(V_2)$ are generated by unitary extensions $W_1$ and $W_2$
(see considerations before the statement of the theorem). Therefore
$\mathbf{R}_{u;z_1}(V_1)$ and $\mathbf{R}_{u;z_2}(V_2)$ are generated by unitary extensions $W_{1,\mathrm{min}}$ and $W_{2,\mathrm{min}}$
as well.
By Lemma~\ref{l1_1}, $(iii)$ we conclude that $S_{j,\mathrm{min}}$ and $W_{j,\mathrm{min}}$ are unitarily equivalent, $j=1,2$.
Thus, $W_{1,\mathrm{min}}$ and $W_{2,\mathrm{min}}$ have no nonzero fixed elements.

Denote 
$$ \mathcal{F}_1 := \{ h\in K\ \arrowvert\ W_1 h = h \}. $$
By the commutativity of $W_1$ and $W_2$, it is easy to check that $\mathcal{F}_1$ reduces $W_2$ (as well as $W_1$). 
Therefore $K\ominus \mathcal{F}_1$ reduce $W_1$ and $W_2$.
Let us check that
\begin{equation}
\label{f3_47}
\mathcal{F}_1 \subseteq K\ominus K_{1,\mathrm{min}}.
\end{equation}
Choose an arbitrary element $h\in \mathcal{F}_1$. Let $h= h_1 + h_2$, where $h_1\in K_{1,\mathrm{min}}$,
$h_2\in K\ominus K_{1,\mathrm{min}}$. Then
$$ h = W_1 h = W_1 h_1 + W_1 h_2, $$
where the first summand on the right belongs to $K_{1,\mathrm{min}}$, while 
the second one belongs to $K\ominus K_{1,\mathrm{min}}$.
Therefore 
$$ h_1 = W_1 h_1,\quad h_2 = W_1 h_2. $$
Since $W_1|_{ K_{1,\mathrm{min}} } = W_{1,\mathrm{min}}$ has no nonzero fixed elements, we get $h_1 = 0$.
So, relation~(\ref{f3_47}) holds. It follows that
\begin{equation}
\label{f3_49}
K_{1,\mathrm{min}} \subseteq K\ominus \mathcal{F}_1,
\end{equation}
and therefore $H\subseteq K\ominus \mathcal{F}_1$. 
Denote $\widehat W_k = W_k|_{ K\ominus \mathcal{F}_1 }$, $k=1,2$.
Unitary operators $\widehat W_1$ and $\widehat W_2$ generate the same generalized resolvent $\mathbf{R}_{u;z_1,z_2} (V_1,V_2)$.
Moreover, $\widehat W_1$ has no nonzero fixed elements.

A similar procedure can be applied for the operators $\widehat W_1$ and $\widehat W_2$, considering a set of all fixed elements for 
$\widehat W_2$:
$$ \mathcal{F}_2 := \{ h\in \widehat K \ \arrowvert\ \widehat W_2 h = h \},\qquad  \widehat K := K\ominus \mathcal{F}_1. $$
Subspaces $\mathcal{F}_2$ and $\widehat K\ominus \mathcal{F}_2$ reduce operators $\widehat W_1$, $\widehat W_2$.
Similarly to~(\ref{f3_47}) we have
\begin{equation}
\label{f3_55}
\mathcal{F}_2 \subseteq \widehat K\ominus K_{2,\mathrm{min}}.
\end{equation}
Then
\begin{equation}
\label{f3_59}
K_{2,\mathrm{min}} \subseteq \widehat K\ominus \mathcal{F}_2,
\end{equation}
and $H\subseteq \widehat K\ominus \mathcal{F}_2$. 
Denote $\widetilde W_k = \widehat W_k|_{ \widehat K\ominus \mathcal{F}_2 }$, $k=1,2$.
Unitary operators $\widetilde W_1$ and $\widetilde W_2$ generate the generalized resolvent $\mathbf{R}_{u;z_1,z_2} (V_1,V_2)$.
Moreover, they have no nonzero fixed elements.
$\Box$

\subsection{A parameterization of all generalized resolvents for commuting symmetric and self-adjoint operators}

We shall first recall a parameterization of all generalized resolvents for commuting isometric and unitary operators.
Taking care of fixed elements of operators and using the results of the previous subsection, we shall obtain a
parameterization of all generalized resolvents in the case of commuting symmetric and self-adjoint operators.

Let $V$ be a closed isometric operator and $U$ be a unitary operator in a Hilbert space $H$, such that
\begin{equation}
\label{f4_10}
V U h = U V h,\quad h\in (U^{-1} D(V))\cap D(V).
\end{equation}

Denote by $\mathcal{S}_{V,U}(\mathbb{D};N_0(V),N_\infty(V))$ a set of all functions $\Phi_{z_1}$ from $\mathcal{S}(\mathbb{D};N_0(V),N_\infty(V))$
which satisfy the following relation:
\begin{equation}
\label{f4_30}
(V\oplus \Phi_{z_1}) U = U (V\oplus \Phi_{z_1}),\qquad z_1\in \mathbb{D}.
\end{equation}

\begin{theorem}
\label{t4_1}(\cite[Theorem 4.1]{cit_15000_Zagorodnyuk_2017_NYJM_Generalized_resolvents})
Let $V$ be a closed isometric operator in a Hilbert space $H$, and $U$ be a unitary operator in $H$. Suppose that
relation~(\ref{f4_10}) holds. 
Then the following statements hold:

\begin{itemize}

\item[(i)] The set of all generalized resolvents of a pair $V,U$ is non-empty if and only if
$\mathcal{S}_{V,U}(\mathbb{D};N_0(V),N_\infty(V))\not= \emptyset$;

\item[(ii)]
Suppose that $\mathcal{S}_{V,U}(\mathbb{D};N_0(V),N_\infty(V))\not= \emptyset$. An arbitrary generalized resolvent
of a pair $V,U$ has the following form:
\begin{equation}
\label{f4_35}
\mathbf{R}_{z_1,z_2} =
( -E_H + 2 \left[ E_H - z_1 (V\oplus \Phi_{z_1}) \right]^{-1} ) U(z_2),\qquad z_1\in\mathbb{D},\ z_2\in \mathbb{T}_e,
\end{equation}
where $\Phi_{z_1}\in \mathcal{S}_{V,U}(\mathbb{D};N_0(V),N_\infty(V))$,
and
\begin{equation}
\label{f4_37}
\mathbf{R}_{z_1,z_2} = \mathbf{R}_{ \frac{1}{\overline{z_1}}, \frac{1}{\overline{z_2}} }^*,\quad 
z_1\in\mathbb{D}_e,\ z_2\in \mathbb{T}_e\backslash\{ 0 \}.
\end{equation}
On the other hand, an arbitrary function $\Phi_{z_1}\in \mathcal{S}_{V,U}(\mathbb{D};N_0(V),N_\infty(V))$ defines
by relations~(\ref{f4_35}),(\ref{f4_37}) a generalized resolvent of a pair $V,U$
(for $z_1\in\mathbb{D}_e$, $z_2=0$ we define $\mathbf{R}_{z_1,z_2}$ by the weak continuity:
$\mathbf{R}_{z_1,0} = w.-\lim\limits_{z_2\to 0} \mathbf{R}_{z_1,z_2}$). 
Moreover, for different operator-valued functions from $\mathcal{S}_{V,U}(\mathbb{D};N_0(V),N_\infty(V))$ there correspond
different generalized resolvents of a pair $V,U$.
\end{itemize}
\end{theorem}

In conditions of Theorem~\ref{t4_1}, if
\begin{equation}
\label{f4_70}
U D(V) = D(V),
\end{equation}
then condition~(\ref{f4_30}) is equivalent to
\begin{equation}
\label{f4_75}
\Phi_{z_1} U g = U \Phi_{z_1} g,\qquad g\in H\ominus D(V),\ z_1\in\mathbb{D}.
\end{equation}
Moreover, in that case the set of all generalized resolvents of $V,U$ is non-empty 
(see~\cite[p. 578]{cit_15000_Zagorodnyuk_2017_NYJM_Generalized_resolvents}).


Let $A_1$ be a closed (not necessarily densely defined) symmetric operator
and $A_2$ be a self-adjoint operator in a Hilbert space $H$, such that
their Cayley's transformations
\begin{equation}
\label{f4_75_1}
V_j := (A_j + i E_H) (A_j - i E_H)^{-1},\qquad j=1,2, 
\end{equation}
satisfy relation~(\ref{f4_10}) with $V=V_1, U=V_2$.
Suppose that there exists a generalized resolvent of operators $A_1,A_2$.
By Proposition~\ref{p2_2} we conclude that the set of generalized resolvents of $V,U$ is non-empty.
Now we can apply Theorem~\ref{t4_1} to $V=V_1, U=V_2$ and obtain a parameterization of all generalized resolvents of $V_1, V_2$.
By Proposition~\ref{p2_2}, in order to obtain generalized resolvents of
$A_1,A_2$ we need to extract those of generalized resolvents of $V,U$,
which are generated by commuting unitary extensions without nonzero fixed elements.

By Theorem~\ref{t3_1} a generalized resolvent $\mathbf{R}_{u;z_1,z_2} (V_1,V_2)$ is generated by~(\ref{f3_15}) 
with some commuting unitary extensions without
nonzero fixed elements if and only if relation~(\ref{f3_35}) holds, where 
$\Phi_1(z_1), \Phi_2(z_2)$ are defined by formulas~(\ref{f3_27}),(\ref{f3_29}).
Since $U=V_2$ is unitary, it has a unique generalized resolvent, and it is generated by $U$. Notice that $\Phi_2(z_2)$ is
the corresponding parameter in Chumakin's formula. Since $U$ has no nozero fixed elements, by Proposition~\ref{p2_1}
we get $\Phi_2(z_2)\in \mathcal{S}_{ \mathrm{adm}; U }(\mathbb{D};N_0(U),N_\infty(U))$.
Thus, relation~(\ref{f3_35}) with $j=2$ holds for all generalized resolvents $\mathbf{R}_{u;z_1,z_2} (V_1,V_2)$.
By~(\ref{f3_27}),(\ref{f4_35}) we can verify that
\begin{equation}
\label{f4_76}
\Phi_{z_1} = \Phi_1(z_1),\qquad z_1\in\mathbb{D}. 
\end{equation}
Consequently, a generalized resolvent $\mathbf{R}_{u;z_1,z_2} (V_1,V_2)$ is generated by~(\ref{f3_15}) 
with some commuting unitary extensions without
nonzero fixed elements if and only if the following relation holds:
\begin{equation}
\label{f4_80}
\Phi_{z_1} \in \mathcal{S}_{ \mathrm{adm}; V_1 }(\mathbb{D};N_0(V_1),N_\infty(V_1)).
\end{equation}

\begin{definition}
\label{d4_1} 
Let $A_1$ be a closed symmetric operator
and $A_2$ be a self-adjoint operator in a Hilbert space $H$, such that
their Cayley's transformations $V_1=V$, $V_2=U$, as in~(\ref{f4_75_1}),
satisfy relation~(\ref{f4_10}).
Denote by $\widetilde{\mathcal{S}}_{V,U}(\mathbb{D};N_0(V),N_\infty(V))$ 
a set of all functions $\Phi_{z_1}$ from $\mathcal{S}_{V,U}(\mathbb{D};N_0(V),N_\infty(V))$
which belong to $\mathcal{S}_{ \mathrm{adm}; V_1 }(\mathbb{D};N_0(V_1),N_\infty(V_1))$.
\end{definition}

\begin{remark}
\label{r4_1}
If $\overline{D(A_1)} = H$, then
$\mathcal{S}_{ \mathrm{adm}; V_1 }(\mathbb{D};N_0(V_1),N_\infty(V_1)) = \mathcal{S}(\mathbb{D};N_0(V_1),N_\infty(V_1))$, 
see Remark~\ref{r2_1}. Since $\Phi_1(z_1)$ is the parameter in Chumakin's formula, we have 
$\Phi_1(z_1)\in\mathcal{S}(\mathbb{D};N_0(V_1),N_\infty(V_1))$.
So, in this case we have 
$$ \widetilde{\mathcal{S}}_{V,U}(\mathbb{D};N_0(V),N_\infty(V)) = \mathcal{S}_{V,U}(\mathbb{D};N_0(V),N_\infty(V)). $$ 
\end{remark}

\begin{theorem}
\label{t4_2} 
Let $A_1$ be a closed (not necessarily densely defined) symmetric operator
and $A_2$ be a self-adjoint operator in a Hilbert space $H$, such that
their Cayley's transformations $V_1=V$, $V_2=U$, as in~(\ref{f4_75_1}),
satisfy relation~(\ref{f4_10}).
Then the following statements hold:

\begin{itemize}

\item[(i)] The set of all generalized resolvents of a pair $A_1,A_2$ is non-empty if and only if
$\widetilde{\mathcal{S}}_{V,U}(\mathbb{D};N_0(V),N_\infty(V))\not= \emptyset$;

\item[(ii)]
Suppose that $\widetilde{\mathcal{S}}_{V,U}(\mathbb{D};N_0(V),N_\infty(V))\not= \emptyset$. An arbitrary generalized resolvent
of a pair $A_1,A_2$ has the following form:
$$ 
\mathbf{R}_{s;\lambda_1,\lambda_2} =
\left( 
E_H - 2 \left[ E_H - \frac{\lambda_1 - i}{\lambda_1 + i} 
\left( V\oplus 
\Phi_{ \frac{\lambda_1 - i}{\lambda_1 + i} }
\right) 
\right]^{-1} 
\right) 
U\left( \frac{\lambda_2 - i}{\lambda_2 + i} 
\right),
$$
\begin{equation}
\label{f4_85}
\lambda_1\in\mathbb{C}_+\backslash\{ i \},\ \lambda_2\in \mathbb{R}_e\backslash\{ \pm i \},
\end{equation}
and
\begin{equation}
\label{f4_95}
\mathbf{R}_{s;\overline{\lambda_1},\overline{\lambda_2}} = \mathbf{R}_{s; \lambda_1, \lambda_2},\qquad 
\lambda_1,\lambda_2\in \mathbb{R}_e,
\end{equation}
where $\Phi_{z_1}\in \widetilde{ \mathcal{S} }_{V,U}(\mathbb{D};N_0(V),N_\infty(V))$.

On the other hand, an arbitrary function $\Phi_{z_1}\in \widetilde{ \mathcal{S} }_{V,U}(\mathbb{D};N_0(V),N_\infty(V))$ defines
by relations~(\ref{f4_85}),(\ref{f4_95}) a generalized resolvent of a pair $A_1,A_2$
(if $\lambda_1=\pm i$, or $\lambda_2=\pm i$, we can define $\mathbf{R}_{s;\lambda_1,\lambda_2}$ by the weak continuity). 
For different operator-valued functions from $\widetilde{ \mathcal{S}_{V,U} }(\mathbb{D};N_0(V),N_\infty(V))$ there correspond
different generalized resolvents of a pair $A_1,A_2$.

\end{itemize}

\end{theorem}

\noindent
\textbf{Proof.} 
$(i)$: Suppose that the set of generalized resolvents of $A_1,A_2$ is non-empty. By~Proposition~\ref{p2_2} for
a generalized resolvent $\mathbf{R}_{s;\lambda_1,\lambda_2}$ of $A_1,A_2$ it corresponds a generalized resolvent
$\mathbf{R}_{u;z_1,z_2}$ of $V_1,V_2$, generated by unitary extensions without non-zero fixed elements.
By considerations before the statement of the theorem this happens if and only if the corresponding
parameter $\Phi_{z_1}$ belongs to $\widetilde{\mathcal{S}}_{V,U}(\mathbb{D};N_0(V),N_\infty(V))$.
Therefore $\widetilde{\mathcal{S}}_{V,U}(\mathbb{D};N_0(V),N_\infty(V))$ is non-empty.

Conversely, suppose that $\widetilde{\mathcal{S}}_{V,U}(\mathbb{D};N_0(V),N_\infty(V))\not= \emptyset$.
Choose an arbitrary $\Phi_{z_1}$ from $\widetilde{\mathcal{S}}_{V,U}(\mathbb{D};N_0(V),N_\infty(V))$.
Then $\Phi_{z_1}\in \mathcal{S}_{V,U}(\mathbb{D};N_0(V),N_\infty(V))$ and
$\Phi_{z_1} \in \mathcal{S}_{ \mathrm{adm}; V_1 }(\mathbb{D};N_0(V_1),N_\infty(V_1))$.
By Theorem~\ref{t4_1} there exists a generalized resolvent $\mathbf{R}_{z_1,z_2} = \mathbf{R}_{u;z_1,z_2}$ of $V,U$ such that
relation~(\ref{f4_35}) holds.
Let us check that $\mathbf{R}_{u;z_1,z_2}$ is generated by some unitary commuting extensions without non-zero fixed elements.
We can repeat the arguments from the proof of sufficiency of Theorem~\ref{t3_1} 
with $\Phi_1(z_1) := \Phi_{z_1}$, $\Phi_2(z_2) := 0$
(we cannot apply this theorem directly). Of course, $\Phi_1(z_1)$, $\Phi_2(z_2)$ satisfy condition~(\ref{f3_35}).
Repeating the construction step-by-step we can get unitary extensions of $V,U$, having no non-zero fixed elements and generating $\mathbf{R}_{u;z_1,z_2}$.
By Proposition~\ref{p2_2} we conclude that the set of generalized resolvents of $A_1,A_2$ is non-empty.

\noindent
$(ii)$: Let $\mathbf{R}_{s;\lambda_1,\lambda_2}$ be an arbitrary generalized resolvent of $A_1,A_2$.
Let $\mathbf{R}_{u;z_1,z_2}$ be the corresponding by Proposition~\ref{p2_2} generalized resolvent, generated
by unitary extensions without fixed elements.
By Proposition~\ref{p2_2} and formula~(\ref{f4_35}) of Theorem~\ref{t4_1} we get relation~(\ref{f4_85})
with some $\Phi_{z_1}\in \mathcal{S}_{V,U}(\mathbb{D};N_0(V),N_\infty(V))$.
Considerations before Definition~\ref{d4_1} show that $\Phi_{z_1}$ satisfy~(\ref{f4_80}).
By Definition~\ref{d4_1} this means that
$\Phi_{z_1}\in \widetilde{ \mathcal{S} }_{V,U}(\mathbb{D};N_0(V),N_\infty(V))$.

Conversely, choose an arbitrary function $\Phi_{z_1}\in \widetilde{ \mathcal{S} }_{V,U}(\mathbb{D};N_0(V),N_\infty(V))$.
Let $\mathbf{R}_{u;z_1,z_2}$ be the corresponding to $\Phi_{z_1}$, by Theorem~\ref{t4_1}, generalized resolvent of $V,U$.
Considerations before Definition~\ref{d4_1} show that $\mathbf{R}_{u;z_1,z_2}$ is generated
by unitary extensions without fixed elements.
By Proposition~\ref{p2_2} $\mathbf{R}_{u;z_1,z_2}$ is related with a generalized resolvent $\mathbf{R}_{s;\lambda_1,\lambda_2}$ of
$A_1,A_2$.
Repeating for $\mathbf{R}_{s;\lambda_1,\lambda_2}$ the construction from the first part of the proof of $(ii)$,
we get formula~(\ref{f4_85}).
Choose an arbitrary $\widehat \Phi_{z_1}\in \widetilde{ \mathcal{S} }_{V,U}(\mathbb{D};N_0(V),N_\infty(V))$,
different from $\Phi_{z_1}$.
By Theorem~\ref{t4_1} the corresponding to $\widehat{ \Phi }_{z_1}$ generalized resolvent $\widehat{ \mathbf{R} }_{u;z_1,z_2}$ 
is different from $\mathbf{R}_{u;z_1,z_2}$. By Proposition~\ref{p2_2} the corresponding generalized resolvents of $A_1,A_2$
will be different as well.
The proof is complete.
$\Box$

\section{The two-dimensional moment problem under Carleman-type conditions}
\label{s_3}

\subsection{The two-dimensional moment problem and generalized resolvents}

We shall need a usual operator-theoretic interpretation of the two-dimensional moment problem, see, e.g.~\cite{cit_980_Fuglede}.
It is rather standard, up to a choice of notations which are due to the taste of a particular researcher.

Let the moment problem~(\ref{f1_1_d_2}) be given with some prescribed set of real numbers $\mathcal{S} := \{ s_{m,n} \}_{m,n=0}^\infty$.
Suppose that this moment problem has a solution $\mu$. Let
\begin{equation}
\label{f5_10}
p(t_1,t_2) = \sum_{m,n=0}^\infty \alpha_{m,n} t_1^m t_2^n,\qquad \alpha_{m,n}\in\mathbb{C},
\end{equation}
where all but finite number of $\alpha_{m,n}$ are zeros.
Calculating $\int |p|^2 d\mu$, one obtains a well-known condition (positive semi-definiteness of $\mathcal{S}$):
\begin{equation}
\label{f5_12}
\sum_{m,n,m',n'=0}^\infty \alpha_{m,n} \overline{ \alpha_{m',n'} } s_{m+m',n+n'}\geq 0,
\end{equation}
for all finitely supported complex sequences $\{ \alpha_{n,m} \}_{m,n=0}^\infty$.

Conversely, suppose that the moment problem~(\ref{f1_1_d_2}) be given with $\mathcal{S} = \{ s_{m,n} \}_{m,n=0}^\infty$, and
$\mathcal{S}$ is positive semi-definite. If $p$ has form~(\ref{f5_10}), 
and $q$ has form~(\ref{f5_10}) with $\beta_{m,n}$ instead of $\alpha_{m,n}$,
then we may consider the following form:
\begin{equation}
\label{f5_15}
S(p,q) = \sum_{m,n,m',n'=0}^\infty \alpha_{m,n} \overline{ \beta_{m',n'} } s_{m+m',n+n'}.
\end{equation}
The form $S$ is sesquilinear and $S(p,p)\geq 0$.
Consider the linear vector space $\mathfrak{L}$ of all polynomials of the form~(\ref{f5_10}), with usual
operations. The linear vector space $\mathfrak{L}$ with $S$ is a pre-Hilbert space.
By the standard procedure, introducing the classes of the equivalence and making the completion, we obtain
a Hilbert space $H$. 
Namely, by $[p] = [p]_{ \mathfrak{L} }$ we denote the set which consists of all $q\in\mathfrak{L}$ such that $S(p-q,p-q) = 0$.
We shall say that $H$ is \textit{the associated Hilbert space to the moment problem~(\ref{f1_1_d_2}) 
(with a positive semi-definite set $\{ s_{m,n} \}_{m,n=0}^\infty$)}.
Denote
\begin{equation}
\label{f5_20}
h_{m,n} = \left[
t_1^m t_2^n
\right]_{ \mathfrak{L} },\qquad  m,n\in\mathbb{Z}_+.
\end{equation}
Observe that
\begin{equation}
\label{f5_25}
(h_{m,n},h_{m',n'})_H = s_{m+m',n+n'},\qquad  m,n,m',n'\in\mathbb{Z}_+.
\end{equation}
If $p$ has form~(\ref{f5_10}), we set
\begin{equation}
\label{f5_27}
A_{1,0} [p]_{ \mathfrak{L} } = \left[
\sum_{m,n=0}^\infty \alpha_{m,n} t_1^{m+1} t_2^n
\right]_{ \mathfrak{L} },\ 
A_{2,0} [p]_{ \mathfrak{L} } = \left[
\sum_{m,n=0}^\infty \alpha_{m,n} t_1^m t_2^{n+1}
\right]_{ \mathfrak{L} }. 
\end{equation}
In particular, we have
\begin{equation}
\label{f5_29}
A_{1,0} h_{m,n} = h_{m+1,n},\quad
A_{2,0} h_{m,n} = h_{m,n+1},\qquad m,n=0,1,2,....
\end{equation}
Let us check that $A_{1,0}$ and $A_{2,0}$ are well-defined. 
If
$q$ has form~(\ref{f5_10}) with $\beta_{m,n}$ instead of $\alpha_{m,n}$, and $q\in [p]$, then
$$ \left(
\sum_{m,n=0}^\infty
\beta_{m,n} h_{m+1,n}, h_{m',n'}
\right)_H
= \sum_{m,n=0}^\infty
\beta_{m,n}
\left(
h_{m+1,n}, h_{m',n'}
\right)_H =
$$
$$
= \sum_{m,n=0}^\infty
\beta_{m,n} s_{m+m'+1,n+n'} =
\sum_{m,n=0}^\infty
\beta_{m,n}
\left(
h_{m,n}, h_{m'+1,n'}
\right)_H = ([p],h_{m'+1,n'}), $$
for arbitrary $m',n'\in\mathbb{Z}_+$.
Since $\{ h_{m',n'} \}_{m',n'=0}^\infty$ are dense in $H$, we see that the value of $A_1$ does not depend on the choice of
a representative in $[p]$. For $A_2$ the argument is similar. Thus, $A_{1,0}$ and $A_{2,0}$ are well-defined.
Moreover, they are linear.
If
$r$ has form~(\ref{f5_10}) with $\gamma_{m,n}$ instead of $\alpha_{m,n}$, then
$$ \left(
A_{1,0} [p], [r]
\right)_H
=
\sum_{m,n,m',n'=0}^\infty \alpha_{m,n} \overline{ \gamma_{m',n'} } s_{m+1+m',n+n'} = \left(
[p], A_{1,0} [r]
\right)_H. $$
Thus, $A_{1,0}$ is symmetric. In a similar manner we obtain that $A_{2,0}$ is symmetric as well.
By definition~(\ref{f5_27}) we see that $A_{1,0}$ and $A_{2,0}$ commute.
Denote $A_1 = \overline{ A_1 }$ and $A_2 = \overline{ A_2 }$.
By~(\ref{f5_29}) and the induction argument we conclude that
\begin{equation}
\label{f5_35}
h_{m,n} = A_1^m A_2^n h_{0,0},\qquad m,n=0,1,2,....
\end{equation}

Suppose now that there exist commuting self-adjoint operators $B_1\supseteq A_1$, $B_2\supseteq A_2$ in a Hilbert space
$\widetilde H\supseteq H$. Then
$$ s_{m,n} = (h_{m,n},h_{0,0})_H = (A_1^m A_2^n h_{0,0},h_{0,0})_H = 
(B_1^m B_2^n h_{0,0},h_{0,0})_{\widetilde H} = $$
$$ = \int_{ \mathbb{R}^2 } t_1^m t_2^n d\widetilde\mu,\qquad  m,n\in\mathbb{Z}_+, $$
where
\begin{equation}
\label{f5_40}
\widetilde\mu(\delta) = (E(\delta) h_{0,0}, h_{0,0})_{\widetilde H},\qquad \delta\in\mathfrak{B}(\mathbb{R}^2),
\end{equation}
and $E(\delta)$ is the operator-valued spectral measure for commuting self-adjoint operators $B_1$ and $B_2$.
Thus, $\widetilde\mu$ is a solution of the moment problem~(\ref{f1_1_d_2}).

It is not hard to see that the above standard construction, which we have recalled, is closely related
to the notion of a generalized resolvent for commuting symmetric operators given in~(\ref{f3_10}).
Since $B_1$ and $B_2$ commute, then their Cayley's transformations commute as well.
Therefore Cayley's transformations $V_j$ of operators $A_j$, $j=1,2$, commute whenever possible.
Thus, if there exists at least one such pair $B_1,B_2$, extending $A_1,A_2$, then condition~(\ref{f3_7}) holds.
In the case, which we are going to study, we shall explicitly construct such a pair $B_1,B_2$.

We may use the generalized resolvent~(\ref{f3_10}), generated by $B_1,B_2$. We can write
$$ \left(
\mathbf{R}_{s;\lambda_1,\lambda_2}(A_1,A_2) h_{0,0}, h_{0,0}
\right)_H = $$
$$ = \left(
(E_{\widetilde H} + \lambda_1 B_1) (B_1 - \lambda_1 E_{\widetilde H})^{-1}
(E_{\widetilde H} + \lambda_2 B_2) (B_2 - \lambda_2 E_{\widetilde H})^{-1}
h_{0,0}, h_{0,0}
\right)_{\widetilde H} = $$
\begin{equation}
\label{f5_45}
= \int_{ \mathbb{R}^2 } 
\frac{ (1+\lambda_1 t_1) (1+\lambda_2 t_2) }{ (t_1 - \lambda_1) (t_2 - \lambda_2) }
d\widetilde\mu,\qquad
\lambda_1, \lambda_2\in \mathbb{R}_e.
\end{equation}
Consequently, \textit{if there exists a generalized resolvent $\mathbf{R}_{s;\lambda_1,\lambda_2}(A_1,A_2)$, then
it generates a solution $\widetilde\mu$ of the moment problem~(\ref{f1_1_d_2}) by the following relation:}
$$ \left(
\mathbf{R}_{s;\lambda_1,\lambda_2}(A_1,A_2) h_{0,0}, h_{0,0}
\right)_H = $$
\begin{equation}
\label{f5_47}
= \int_{ \mathbb{R}^2 } 
\frac{ (1+\lambda_1 t_1) (1+\lambda_2 t_2) }{ (t_1 - \lambda_1) (t_2 - \lambda_2) }
d\widetilde\mu(t_1,t_2),\qquad
\lambda_1, \lambda_2\in \mathbb{R}_e.
\end{equation}
Let us explain a possible way to get $\widetilde\mu$ from relation~(\ref{f5_47}).
Consider the following transformation $T$, which maps $\mathbb{R}^2$ onto $(0,2\pi]^2$:
\begin{equation}
\label{f5_49}
\mathbb{R}^2 \ni (t_1,t_2) \mapsto T ((t_1,t_2)) = (\theta_1,\theta_2)\in (0,2\pi]^2,\
\frac{ t_k + i }{ t_k - i } = e^{ i\theta_k },\ k=1,2,
\end{equation}
and 
\begin{equation}
\label{f5_51}
\sigma(\delta) := \widetilde\mu (T^{-1} \delta),\qquad \delta\in\mathfrak{B}((0,2\pi]^2),
\end{equation}
where $T^{-1}$ is the inverse of $T$. Since
$$ \frac{ 1+ \lambda_k t_k }{ t_k - \lambda_k } = 
i \frac{ 1 + z_k e^{i\theta_k} }{  1 - z_k e^{i\theta_k} },\quad k=1,2, $$
then the right-hand side of~(\ref{f5_47}) is equal to
\begin{equation}
\label{f5_55}
- \int_{ (0,2\pi]^2 } 
\frac{ (1 + z_1 e^{i\theta_1}) (1 + z_2 e^{i\theta_2}) }{  (1 - z_1 e^{i\theta_1}) (1 - z_2 e^{i\theta_2}) } 
d\sigma(\theta_1,\theta_2),\qquad
z_1, z_2\in \mathbb{T}_e.
\end{equation}
The latter function, by taking derivatives, uniquely determines all trigonometric moments of the Borel measure $\sigma$,
and therefore determines $\sigma$, see, e.g.~\cite{cit_15000_Zagorodnyuk_2017_NYJM_Generalized_resolvents}.
Then $\widetilde\mu$ is given by relation~(\ref{f5_51}).

Denote by $\tau$ the correspondence which relates a solution $\widetilde\mu$ to an arbitrary generalized resolvent
$\mathbf{R}_{s;\lambda_1,\lambda_2}(A_1,A_2)$ by relation~(\ref{f5_47}).

\begin{proposition}
\label{p5_1}
Let the moment problem~(\ref{f1_1_d_2}) be given with some prescribed set of real numbers $\mathcal{S} = \{ s_{m,n} \}_{m,n=0}^\infty$,
and condition~(\ref{f5_12}) holds for all finitely supported complex sequences $\{ \alpha_{m,n} \}_{m,n=0}^\infty$.
Let $H$ be the associated Hilbert space and $A_1,A_2$ be the corresponding symmetric operators in $H$.
Then the following statements hold:
\begin{itemize}
\item[(i)] The moment problem~(\ref{f1_1_d_2}) has a solution if and only if there exists a generalized resolvent 
$\mathbf{R}_{s;\lambda_1,\lambda_2}(A_1,A_2)$;

\item[(ii)] If the moment problem~(\ref{f1_1_d_2}) is solvable, then 
there exists a $1-1$ correspondence between a set of all generalized resolvents $\mathbf{R}_{s;\lambda_1,\lambda_2}(A_1,A_2)$
and a set of all all solutions $\mu$ of the moment problem~(\ref{f1_1_d_2}).
The latter correspondence is provided by $\tau$.
\end{itemize}

\end{proposition}

\noindent
\textbf{Proof.} 
Suppose that there exists a solution $\mu$ of the moment problem~(\ref{f1_1_d_2}). 
Consider the usual space $L^2_\mu$ of (the classes of the equivalence of) square summable complex-valued functions on $\mathbb{R}^2$.
The class of the equivalence, which contains $f(t_1,t_2)$, will be denoted by $[f] = [f]_{L^2_\mu}$.
Let $\mathcal{B}_k$ be the operator in $L^2_\mu$ of multiplication of $f(t_1,t_2)$ by $t_k$, $k=1,2$.
Operators $\mathcal{B}_1,\mathcal{B}_2$ are commuting self-adjoint operators.
Denote by $P^2_{\mu}$ a subset of $L^2_\mu$ which consists of those classes of equivalence which contain polynomials,
and $L^2_{\mu,0} := \overline{ P^2_{\mu} }$.
Consider the following transformation $W_0$ of $P^2_{\mu}$ into $H$:
\begin{equation}
\label{f5_57}
W_0 \left[
p(t_1,t_2)
\right]_{L^2_\mu}
=
\sum_{m,n=0}^\infty \alpha_{m,n} h_{m,n},
\end{equation}
where $p$ is from~(\ref{f5_10}).
Let us check that this transformation is well-defined.
Suppose that the class $[p]_{L^2_\mu}$ contains another polynomial $q(t_1,t_2)$, which has form~(\ref{f5_10}), with some complex
numbers $\beta_{m,n}$ instead of $\alpha_{m,n}$. Then

$$ \left\|
\sum_{m,n=0}^\infty \alpha_{m,n} h_{m,n} - \sum_{m,n=0}^\infty \beta_{m,n} h_{m,n} 
\right\|_H^2 = $$
$$ = \left(
\sum_{m,n=0}^\infty (\alpha_{m,n} - \beta_{m,n}) h_{m,n}, \sum_{m',n'=0}^\infty (\alpha_{m',n'} - \beta_{m',n'}) h_{m',n'} 
\right)_H = $$
$$ = 
\sum_{m,n,n',n'=0}^\infty (\alpha_{m,n} - \beta_{m,n}) \overline{ (\alpha_{m',n'} - \beta_{m',n'}) }
s_{m+m',n+n'} = $$
$$ = \left\|
[ p(t_1,t_2) ]_{L^2_\mu} - [ q(t_1,t_2) ]_{L^2_\mu}
\right\|_{L^2_\mu}^2 = 0. $$
Thus, $W_0$ is well-defined. 
Moreover, $W_0$ is linear. 
With $p$ from~(\ref{f5_10}) and $r$ of the same form, with some complex $\gamma_{m,n}$ instead of
$\alpha_{m,n}$, we may write
$$ (W_0 [p], W_0 [r])_H = 
\left(
\sum_{m,n=0}^\infty \alpha_{m,n} h_{m,n}, \sum_{m',n'=0}^\infty \gamma_{m,n} h_{m,n}
\right)_H = $$
$$ = 
\sum_{m,n,n',n'=0}^\infty \alpha_{m,n} \overline{ \gamma_{m',n'} }
s_{m+m',n+n'} = ([p],[r])_{ L^2_\mu }. $$
Therefore $W_0$ is isometric. We extend $W_0$ by the continuity to a unitary transformation $W$, which maps
$L^2_{\mu,0}$ on the whole space $H$.
Set
\begin{equation}
\label{f5_59}
U := W \oplus E_{ L^2_\mu\ominus L^2_{\mu,0} }.
\end{equation}
Observe that $U$ is a unitary transformation which maps $L^2_\mu$ onto $\widetilde H := H\oplus \left( L^2_\mu\ominus L^2_{\mu,0} \right)$.
Denote
\begin{equation}
\label{f5_62}
\widetilde A_1 = U \mathcal{B}_1 U^{-1},\quad \widetilde A_2 = U \mathcal{B}_2 U^{-1}.
\end{equation}
Operators $\widetilde A_1$ and $\widetilde A_2$ are commuting self-adjoint operators in $\widetilde H\supseteq H$.
Since
$$ \widetilde A_1 h_{m,n} = U \mathcal{B}_1 [t_1^m t_2^n] = U [t_1^{m+1} t_2^n] = h_{m+1,n}, $$
$$ \widetilde A_2 h_{m,n} = U \mathcal{B}_2 [t_1^m t_2^n] = U [t_1^m t_2^{n+1}] = h_{m,n+1},\qquad m,n\in\mathbb{Z}_+, $$
then $\widetilde A_k\supseteq A_k$, $k=1,2$.
Denote by $\mathbf{R}_{s;\lambda_1,\lambda_2}$ the generalized resolvent of $A_1,A_2$, which is generated by 
extensions $\widetilde A_1, \widetilde A_2$.
The existence of $\mathbf{R}_{s;\lambda_1,\lambda_2}$ proves the necessity of statement~$(i)$, while its sufficiency follows from
considerations before the statement of the proposition. Now we proceed to check statement~$(ii)$.
By relation~(\ref{f5_47}), $\mathbf{R}_{s;\lambda_1,\lambda_2}$ produces a solution $\widetilde\mu$ of 
the moment problem~(\ref{f1_1_d_2}). Let us check that $\widetilde\mu = \mu$.
The left-hand side of~(\ref{f5_47}) is equal to
$$ \left(
(E_{\widetilde H} + \lambda_1 \widetilde A_1) (\widetilde A_1 - \lambda_1 E_{\widetilde H})^{-1}
(E_{\widetilde H} + \lambda_2 \widetilde A_2) (\widetilde A_2 - \lambda_2 E_{\widetilde H})^{-1}
h_{0,0}, h_{0,0}
\right)_{\widetilde H} = $$
$$ = 
\left(
(E_{L^2_\mu} + \lambda_1 \mathcal{B}_1) (\mathcal{B}_1 - \lambda_1 E_{L^2_\mu})^{-1}
(E_{L^2_\mu} + \lambda_2 \mathcal{B}_2) (\mathcal{B}_2 - \lambda_2 E_{L^2_\mu})^{-1}
[1], [1]
\right)_{L^2_\mu} = $$
\begin{equation}
\label{f5_64}
= \int_{ \mathbb{R}^2 } 
\frac{ (1+\lambda_1 t_1) (1+\lambda_2 t_2) }{ (t_1 - \lambda_1) (t_2 - \lambda_2) }
d(F(\delta) [1], [1])_{L^2_\mu},\qquad
\lambda_1, \lambda_2\in \mathbb{R}_e,
\end{equation}
where 
$F(\delta)$ is the spectral measure for commuting self-adjoint operators $\mathcal{B}_1$ and $\mathcal{B}_2$.
Since $F(\delta)$ is the operator of multiplication by the characteristic function $\chi(\delta)$ in $L^2_\mu$
(notice that $L^2_\mu$ is a particular case of the direct integral of Hilbert spaces and multiplication operators have the required property therein),
then
$$ (F(\delta) [1], [1])_{L^2_\mu} = \mu(\delta),\qquad \delta\in\mathfrak{B}(\mathbb{R}^2). $$
As it was discussed above, the integrals of the form~(\ref{f5_64}) determine the measure uniquely.
Consequently, we have $\mu = \widetilde\mu$.
Thus, $\tau$ is surjective.

Suppose that two generalized resolvents $\mathbf{R}^{[k]}_{s;\lambda_1,\lambda_2}$ 
of $A_1,A_2$, $k=1,2$, generate 
by relation~(\ref{f5_47}) the same solution $\widetilde\mu$.
Suppose that $\mathbf{R}^{[k]}_{s;\lambda_1,\lambda_2}$  is generated by commuting
self-adjoint operators $B_1^{[k]}$ and $B_2^{[k]}$ in a Hilbert space $\widetilde H^{[k]}$, $k=1,2$.
For arbitrary $m,m',n,n'\in\mathbb{Z}_+$, and $k=1,2$, we may write
$$ \left(
\mathbf{R}^{[k]}_{s;\lambda_1,\lambda_2} h_{m,n}, h_{m',n'}
\right)_H
=
\left(
\mathbf{R}^{[k]}_{s;\lambda_1,\lambda_2} A_1^m A_2^n h_{0,0}, A_1^{m'} A_2^{n'} h_{0,0}
\right)_H
=
$$
$$
= \left(
(E + \lambda_1 B^{[k]}_1) (B^{[k]}_1 - \lambda_1 E)^{-1}
(E + \lambda_2 B^{[k]}_2) (B^{[k]}_2 - \lambda_2 E)^{-1} * \right. $$
$$ \left. *
(B^{[k]}_1)^{m+m'} (B^{[k]}_2)^{n+n'}
h_{0,0}, h_{0,0}
\right)_{\widetilde H^{[k]}} = $$
$$ = \int_{ \mathbb{R}^2 } 
\frac{ (1+\lambda_1 t_1) (1+\lambda_2 t_2) }{ (t_1 - \lambda_1) (t_2 - \lambda_2) } t_1^{m+m'} t_2^{n+n'}
d\widetilde\mu(t_1,t_2),\qquad
\lambda_1, \lambda_2\in \mathbb{R}_e.
$$
Here we have used relation~(\ref{f5_40}) for the measure $\widetilde\mu$, which ensures that
$\widetilde\mu$ is generated by the spectral measures of extensions.
Thus, we have 
$$ \mathbf{R}^{[1]}_{s;\lambda_1,\lambda_2} = \mathbf{R}^{[2]}_{s;\lambda_1,\lambda_2}. $$
So, $\tau$ is injective as well.
The proof is now complete. $\Box$

\subsection{The two-dimensional moment problem when the associated operator $A_2$ is self-adjoint}

Let the moment problem~(\ref{f1_1_d_2}) be given with a prescribed set of real numbers $\mathcal{S} = \{ s_{m,n} \}_{m,n=0}^\infty$.
and condition~(\ref{f5_12}) holds for all finitely supported complex sequences $\{ \alpha_{n,m} \}_{m,n=0}^\infty$.
Let $H$ be the associated Hilbert space and $A_1,A_2$ be the corresponding symmetric operators in $H$.
We shall need conditions which ensure that the operator $A_2$ is self-adjoint.
There are different effective methods to determine that a given operator is self-adjoint, see, e.g.,~\cite{cit_1000_Berezansky_1965__Book}.
One of relevant tools is the use of quasi-analiticity, which leads to Carleman-type conditions.

Let $A$ be a symmetric operator in a Hilbert space $H$. Denote
$$ D^\infty(A) = \bigcap_{n=1}^\infty D(A^n),\quad D_{\mathrm{qa}}(A) = 
\{
x\in D^\infty(A)\ \arrowvert\ \sum_{n=0}^\infty \| A^n x \|^{-1/n} = \infty
\}. $$
Vectors from $D_{\mathrm{qa}}(A)$ are called \textit{quasi-analytic}. The latter notion was introduced by Nussbaum who
proved the following theorem:

\begin{theorem}
\label{t6_1}
(\cite[Theorem QA]{cit_983__Chernoff_1972})
Let $A$ be a symmetric operator on Hilbert space. If $A$ has a total set of quasi-analytic vectors then the closure of $A$ is
selfadjoint.
\end{theorem}
Moreover, if $x\in D_{\mathrm{qa}}(A)$, then $A^k x\in D_{\mathrm{qa}}(A)$, $\forall k\in\mathbb{N}$ (see
Lemma~2.1 in~\cite{cit_983__Chernoff_1972}).

Let us check that condition~(\ref{f1_2}) implies self-adjointness of $A_2$. In fact, by Theorem~\ref{t6_1} and the remark after it,
$A_2$ is self-adjoint if vectors
$\{ h_{m,0} \}_{m=0}^\infty$ are quasi-analytic.
Vector $h_{m,0}$ is quasi-analytic if
$$ \infty = \sum_{k=0}^\infty \| A_2^k h_{m,0} \|^{-1/k} = \sum_{k=0}^\infty (h_{m,k}, h_{m,k})^{-1/(2k)} =
\sum_{k=0}^\infty \frac{1}{ \sqrt[2k]{s_{2m,2k}} }. $$
By the comparison test we see that condition~(\ref{f1_2}) implies that the latter relation holds
and therefore $A_2$ is self-adjoint.
We see that the following condition:
\begin{equation}
\label{f6_5}
\sum_{k=0}^\infty \frac{1}{ \sqrt[2k]{s_{2m,2k}} } = \infty,\quad \forall m\in\mathbb{Z}_+,
\end{equation}
implies self-adjointness of $A_2$. However, we shall need condition~(\ref{f1_2}) to ensure the existence of
generalized resolvents for $A_1,A_2$.

Besides operators $A_1$ and $A_2$ in $H$, we shall use the following operator $J$.
If $p$ has form~(\ref{f5_10}), we set
\begin{equation}
\label{f65_10}
J_0 [p]_{ \mathfrak{L} } = \left[
\sum_{m,n=0}^\infty \overline{ \alpha_{m,n} } t_1^m t_2^n
\right]_{ \mathfrak{L} }. 
\end{equation}
In particular, we have
\begin{equation}
\label{f6_15}
J_{0} h_{m,n} = h_{m,n},\qquad m,n=0,1,2,....
\end{equation}

Let us check that $J_0$  is well-defined. If
$q$ has form~(\ref{f5_10}) with some complex $\beta_{m,n}$ instead of $\alpha_{m,n}$, and $q\in [p]$, then

$$ \left\|
\left[
\sum_{m,n=0}^\infty \overline{ \alpha_{m,n} } t_1^m t_2^n
\right]_{ \mathfrak{L} }
-
\left[
\sum_{m,n=0}^\infty \overline{ \beta_{m,n} } t_1^m t_2^n
\right]_{ \mathfrak{L} }
\right\|_H^2 = $$
$$ \left(
\left[
\sum_{m,n=0}^\infty \overline{ (\alpha_{m,n} - \beta_{m,n}) }  t_1^m t_2^n
\right]_{ \mathfrak{L} },
\left[
\sum_{m',n'=0}^\infty \overline{ (\alpha_{m',n'} - \beta_{m',n'}) } t_1^{m'} t_2^{n'}
\right]_{ \mathfrak{L} }
\right)_H^2 = $$
$$ = 
\sum_{m,n,m',n'=0}^\infty \overline{ (\alpha_{m,n} - \beta_{m,n}) } (\alpha_{m',n'} - \beta_{m',n'})
s_{m+m',n+n'} = $$
$$ = 
\sum_{m,n,m',n'=0}^\infty \overline{ (\alpha_{m,n} - \beta_{m,n}) } (\alpha_{m',n'} - \beta_{m',n'})
\left(
\left[
t_1^{m'} t_2^{n'} \right]_{ \mathfrak{L} }, 
\left[ 
t_1^m t_2^n
\right]_{ \mathfrak{L} }
\right)_H = $$
$$ = \left\|
[p(t_1,t_2)] - [q(t_1,t_2)]
\right\|_H^2 = 0. $$
Thus, $J_0$ is well-defined. For $p$ as in~(\ref{f5_10}), and $r$ of the same form, with some complex $\gamma_{m,n}$ instead of
$\alpha_{m,n}$, we may write
$$ (J_0 [p], J_0 [r])_H = 
\left(
\left[
\sum_{m,n=0}^\infty \overline{ \alpha_{m,n} } t_1^m t_2^n
\right]_{ \mathfrak{L} }, 
\left[
\sum_{m',n'=0}^\infty \overline{ \gamma_{m',n'} } t_1^{m'} t_2^{n'}
\right]_{ \mathfrak{L} }
\right)_H = $$
$$ = 
\sum_{m,n,m',n'=0}^\infty \overline{ \alpha_{m,n} } \gamma_{m',n'} 
s_{m+m',n+n'} = ([r],[p])_H. $$
Notice that $J_0$ is antilinear and $J_0^2 = E_H$. Thus, we may extend $J_0$ by the continuity to a conjugation $J$ on the whole space $H$.

In what follows we shall assume that condition~(\ref{f1_2}) holds true. Let us check that all conditions of
the following lemma are satisfied with
$$ \mathbf{A} = A_{1,0},\ \mathbf{B} = A_{2,0},\ \mathbf{J} = J,\ \mathbf{H} = H,\  \mathcal{D} = \mathfrak{L},\ z_0=i. $$

\begin{lemma}
\label{l6_1}
(\cite{cit_14000_Zagorodnyuk_2012_MFAT__commuting_operators})
Let $\mathbf{A}$ be a symmetric operator and $\mathbf{B}$ be an essentially self-adjoint operator with a common domain
$\mathcal{D} = D(\mathbf{A}) = D(\mathbf{B})$ in a Hilbert space $\mathbf{H}$, $\overline{ \mathcal{D} } = \mathbf{H}$, and
$$ \mathbf{A} \mathcal{D} \subseteq \mathcal{D},\ \mathbf{B} \mathcal{D} \subseteq \mathcal{D}; $$ 
$$ \mathbf{A} \mathbf{B} = \mathbf{B} \mathbf{A}. $$
Suppose also that for some $z_0\in\mathbb{C}\backslash\mathbb{R}$, the operator $\mathbf{B}$ restricted to the domain
$(\mathbf{A} - z_0 E_{\mathbf{H}}) \mathcal{D}$ is essentially self-adjoint in a Hilbert space 
$(\overline{ \mathbf{A} } - z_0 E_{\mathbf{H}}) D(\overline{ \mathbf{A} })$.

If there exists a conjugation $\mathbf{J}$ in $\mathbf{H}$ such that $\mathbf{J} \mathcal{D} \subseteq \mathcal{D}$,
and
$$ \mathbf{A} \mathbf{J} = \mathbf{J} \mathbf{A},\ \mathbf{B} \mathbf{J} = \mathbf{J} \mathbf{B}, $$
then there exists a self-adjoint operator $\widetilde{\mathbf{A}}\supseteq\mathbf{A}$, which commutes with $\overline{ \mathbf{B} }$.
\end{lemma}
Denote
\begin{equation}
\label{f6_20}
\mathbf{H}_1 = (A_1 - i E_H) D(A_1),\ \mathbf{H}_2 = \mathbf{H}\ominus\mathbf{H}_1,\ 
\mathbf{H}_3 = (A_1 + i E_H) D(A_1),\ \mathbf{H}_4 = \mathbf{H}\ominus\mathbf{H}_3.
\end{equation}
Let us check that $A_{2,0}|_{ (A_{1,0} - i E_H) \mathfrak{L}  }$ is essentially self-adjoint in $\mathbf{H}_1$.
By Theorem~\ref{t6_1} and a remark after it, the latter holds if vectors
$$ \{ (A_1 - iE_H) h_{m,0} \}_{m=0}^\infty = \{ h_{m+1,0} - i h_{m,0} \}_{m=0}^\infty $$
are quasianalytic in $\mathbf{H}_1$. The latter means that
\begin{equation}
\label{f6_25}
\sum_{k=0}^\infty \left\|
A_2^k (h_{m+1,0} - i h_{m,0})
\right\|^{-1/k} = \infty,\qquad \forall m\in\mathbb{Z}_+.
\end{equation}
The summand in~(\ref{f6_25}) is equal to
$$ \left\| h_{m+1,k} - i h_{m,k} \right\|^{-1/k} = (h_{m+1,k} - i h_{m,k}, h_{m+1,k} - i h_{m,k})^{-1/(2k)} = $$
$$ = \frac{1}{ \sqrt[2k]{ s_{2m,2k} + s_{2m+2,2k} } }. $$
By~(\ref{f1_2}) we conclude that relation~(\ref{f6_25}) holds. 
Thus, $A_{2,0}|_{ (A_{1,0} - i E_H) \mathfrak{L}  }$ is essentially self-adjoint.
The rest of conditions of Lemma~\ref{l6_1} are easy to verify.

We shall use some constructions from the proof of Lemma~\ref{l6_1}, as given in~\cite{cit_14000_Zagorodnyuk_2012_MFAT__commuting_operators}.
Let $V_{\mathbf{A}} = V = V_1, U_{\mathbf{B}} = U=V_2$ be Cayley's transformations of $A_1,A_2$ from~(\ref{f3_5})
(the notation $V_{\mathbf{A}},U_{\mathbf{B}}$ was used in~\cite{cit_14000_Zagorodnyuk_2012_MFAT__commuting_operators}).
In the above-mentioned proof in~\cite{cit_14000_Zagorodnyuk_2012_MFAT__commuting_operators} it was shown that the unitary operator 
$U_{\mathbf{B}} = V_2$
reduces subspaces $\mathbf{H}_j$, $j=1,2,3,4$. Denote by $W_j = U_{\mathbf{B},j}$ the restriction of $U_{\mathbf{B}}$ to 
$\mathbf{H}_j$, which is a unitary
operator in $\mathbf{H}_j$, $j=1,2,3,4$.
By the Godich-Lutsenko theorem, there exist two conjugations $K,L$ in $\mathbf{H}_2$, such that $W_2 = KL$. Set
$U_{2,4} = \mathbf{J} K$, and $\widetilde U = V_{\mathbf{A}}\oplus U_{2,4}$. Let 
$\widetilde{\mathbf{A}}$ be the inverse Cayley's transformation of $\widetilde U$.
The operator $\widetilde{\mathbf{A}}$ is the operator provided by Lemma~\ref{l6_1}.

Thus, Lemma~\ref{l6_1} guarantees the existence of self-adjoint commuting extensions of $A_1,A_2$ inside $H$. 
As we already
noticed, this means that condition~(\ref{f3_7}) holds, and we may consider the corresponding generalized resolvent.
By Proposition~\ref{p5_1} the moment problem~(\ref{f1_1_d_2}) has a solution.
In what follows we are going to construct all solutions of the moment problem~(\ref{f1_1_d_2}).
At first, we shall describe canonical solutions.


\begin{definition}
\label{d6_1} 
A solution of the moment problem~(\ref{f1_1_d_2}) is said to be \textbf{canonical} if it is generated
by commuting self-adjoint extensions inside the associated Hilbert space $H$ for the corresponding symmetric operators $A_1,A_2$.
\end{definition}

In~\cite{cit_14000_Zagorodnyuk_2012_MFAT__commuting_operators} it was shown that all unitary operators $\mathcal{U}\supseteq V_1$,
commuting with $V_2$, have the following form:
\begin{equation}
\label{f6_30}
\mathcal{U} = V_1 \oplus \left( U_{2,4} U_2 \right),
\end{equation}
where $U_2$ is an arbitrary unitary operator in $\mathbf{H}_2$, which commutes with $V_2$. Such operators $U_2$
are represented as unitary decomposable operators in the direct integral model space for $V_2$.

\begin{theorem}
\label{t6_2}
Let the moment problem~(\ref{f1_1_d_2}) be given with some prescribed set of real numbers $\{ s_{m,n} \}_{m,n=0}^\infty$.
Suppose that conditions~(\ref{f5_12}) and~(\ref{f1_2}) hold.
Let $H$ be the associated Hilbert space, $A_1,A_2$ be the corresponding symmetric operators in $H$ and
$V=V_1,U=V_2$ be their Cayley's transformations from~(\ref{f3_7}).
An arbitrary canonical solution $\mu$ of the moment problem~(\ref{f1_1_d_2}) satisfy the following relation:
$$ \left(
(E_{H} + \lambda_1 \widetilde A_1) (\widetilde A_1 - \lambda_1 E_{H})^{-1}
(E_{H} + \lambda_2 A_2) (A_2 - \lambda_2 E_{H})^{-1}
h_{0,0}, h_{0,0}
\right)_{H} = $$
\begin{equation}
\label{f6_32}
= \int_{ \mathbb{R}^2 } 
\frac{ (1+\lambda_1 t_1) (1+\lambda_2 t_2) }{ (t_1 - \lambda_1) (t_2 - \lambda_2) }
d\mu(t_1,t_2),\qquad
\lambda_1, \lambda_2\in \mathbb{R}_e,
\end{equation}
where
\begin{equation}
\label{f6_34}
\widetilde A_1 = i E_H + 2i \left( 
V_1 \oplus ( U_{2,4} U_2 ) -E_H
\right)^{-1}
\end{equation}
and $U_2$ is an arbitrary unitary operator in $\mathbf{H}_2$, which commutes with $V_2$.

On the other hand, each unitary operator $U_2$ in $\mathbf{H}_2$, which commutes with $V_2$,
generates by~(\ref{f6_34}) and (\ref{f6_32}) a canonical solution $\mu$ of the moment problem.
Moreover, different operators $U_2$ generate different solutions of the moment problem.
\end{theorem}

\noindent
\textbf{Proof. }
All statements of the theorem, except the last one, follow easily from considerations before the statement of the theorem and
Proposition~\ref{p5_1}.
Let us check the last statement.
Since relation~(\ref{f6_34}) uniquely determines $U_2$ from the knowledge of $\widetilde A_!$, it follows that
different $U_2$ can not produce the same operator $\widetilde A_1$.
Denote by $\mathbf{R}_{s;\lambda_1,\lambda_2}$ the generalized resolvent of $A_1,A_2$, generated by $\widetilde A_1, A_2$.
Since
$$ \mathbf{R}_{s;2i,i} = -2 E_H + (-3i)  (\widetilde A_1 - 2i E_H)^{-1}, $$
then $\widetilde A_1$ is uniquely determined by the corresponding generalized resolvent.
Therefore different self-adjoint operators $\widetilde A_1$ can not produce the same generalized resolvent of $A_1,A_2$.
Finally, by Proposition~\ref{p5_1} different generalized resolvents give different solutions of the moment problem.
The proof is now complete. $\Box$


We are going to use Theorem~\ref{t4_2} which provides a parameterization of all generalized resolvents of $A_1,A_2$.
In our case we have some simplifications.
Since $D(A_1) = \mathfrak{L}$ is dense in $H$, then by Remark~\ref{r4_1} we conclude that
$$ \widetilde{\mathcal{S}}_{V,U}(\mathbb{D};N_0(V),N_\infty(V)) = \mathcal{S}_{V,U}(\mathbb{D};N_0(V),N_\infty(V)). $$ 
Since the restriction of $V_2$ to $\mathbf{H}_1$ is a unitary operator, then
$V_2 \mathbf{H}_1 = \mathbf{H}_1$. Since $D(V_1) = \mathbf{H}_1$, we see that condition~(\ref{f4_70}) holds.
As it was noticed below formula~(\ref{f4_70}), condition~(\ref{f4_30}) is then equivalent to the following condition,
which we rewrite here for convenience of reading:
\begin{equation}
\label{f6_50}
\Phi_{z_1} U g = U \Phi_{z_1} g,\qquad g\in H\ominus D(V),\ z_1\in\mathbb{D}.
\end{equation}
In our case, the set $\mathcal{S}_{V,U}(\mathbb{D};N_0(V),N_\infty(V))$ coincides with a set of all functions 
$\Phi_{z_1}$ from $\mathcal{S}(\mathbb{D};N_0(V),N_\infty(V))$
which satisfy relation~(\ref{f6_50}).
Observe that $\Phi_{z_1}$ belongs to $\mathcal{S}_{V,U}(\mathbb{D};N_0(V),N_\infty(V))$ if and only if
$G_{z_1} := \Phi_{z_1} U_{2,4}^{-1}$ belongs to $\mathcal{S}(\mathbb{D};N_\infty(V),N_\infty(V))$, and
$G_{z_1}$ commutes with $U|_{N_\infty(V)}$, for all $z_1\in\mathbb{D}$.
Thus, the values of $G_{z_1}$ are represented as contractive decomposable operators in the direct integral model space for $U|_{N_\infty(V)}$.

\begin{theorem}
\label{t6_3}
Let the moment problem~(\ref{f1_1_d_2}) be given with some set of real numbers $\{ s_{m,n} \}_{m,n=0}^\infty$.
Suppose that conditions~(\ref{f5_12}) and~(\ref{f1_2}) hold.
Let $H$ be the associated Hilbert space, $A_1,A_2$ be the corresponding symmetric operators in $H$ and
$V=V_1,U=V_2$ be their Cayley's transformations from~(\ref{f3_7}).
An arbitrary solution $\mu$ of the moment problem~(\ref{f1_1_d_2}) satisfy the following relation:
$$ \left(
\left( 
E_H - 2 \left[ E_H - \frac{\lambda_1 - i}{\lambda_1 + i} (V\oplus \Phi_{ \frac{\lambda_1 - i}{\lambda_1 + i} }) \right]^{-1} 
\right) 
U\left( \frac{\lambda_2 - i}{\lambda_2 + i} 
\right)
h_{0,0}, h_{0,0}
\right)_{H} = $$
\begin{equation}
\label{f6_52}
= \int_{ \mathbb{R}^2 } 
\frac{ (1+\lambda_1 t_1) (1+\lambda_2 t_2) }{ (t_1 - \lambda_1) (t_2 - \lambda_2) }
d\mu(t_1,t_2),\qquad
\lambda_1\in\mathbb{C}_+\backslash\{ i \},\ \lambda_2\in \mathbb{R}_e\backslash\{ \pm i \},
\end{equation}
where $\Phi_{z_1}\in \mathcal{S}_{V,U}(\mathbb{D};N_0(V),N_\infty(V))$.

On the other hand, each $\Phi_{z_1}\in \mathcal{S}_{V,U}(\mathbb{D};N_0(V),N_\infty(V))$
generates by~(\ref{f6_52}) a solution $\mu$ of the moment problem~(\ref{f1_1_d_2}).
Moreover, different functions $\Phi_{z_1}$ from $\mathcal{S}_{V,U}(\mathbb{D};N_0(V),N_\infty(V))$ generate different solutions of 
the moment problem~(\ref{f1_1_d_2}).
\end{theorem}

\begin{remark}
\label{6_1}
In order to obtain all values of the integral in~(\ref{f6_52}), for $\lambda_1,\lambda_2\in\mathbb{R}_e$
(which one needs to recover $\mu$), one can calculate the complex conjugate values and use the continuity
of the left-hand side. 
\end{remark}

\noindent
\textbf{Proof.} All statements of the theorem follow from Proposition~\ref{p5_1} and
a description of all generalized resolvents of $A_1,A_2$ from Theorem~\ref{t4_2}.
$\Box$

\begin{corollary}
\label{c6_1}
In conditions of Theorem~\ref{t6_3} the following statements are equivalent:

\begin{itemize}
\item[(i)] The moment problem~(\ref{f1_1_d_2}) is determinate;

\item[(ii)] There exists a unique canonical solution of he moment problem~(\ref{f1_1_d_2});

\item[(iii)] The operator $A_1$ is self-adjoint.
\end{itemize}

\end{corollary}

\noindent
\textbf{Proof.} 
$(i)\Rightarrow (ii)$: This is obvious.

\noindent
$(ii)\Rightarrow (iii)$: Suppose that $A_1$ is not self-adjoint. Then $\mathbf{H}_2\not= \{ 0 \}$. Operators $E_{ \mathbf{H}_2 }$ and
$i E_{ \mathbf{H}_2 }$ are different unitary operators in $\mathbf{H}_2$, which generate by Theorem~\ref{t6_2} different
canonical solutions of the moment problem. This contradiction proves the assertion.

\noindent
$(iii)\Rightarrow (i)$: If $A_1$ is self-adjoint, then $\mathbf{H}_4 = \{ 0 \}$, $\mathbf{H}_4 = \{ 0 \}$.
In this case the set $\mathcal{S}(\mathbb{D};N_0(V),N_\infty(V))$ consists of a unique zero function.
Then $\mathcal{S}_{V,U}(\mathbb{D};N_0(V),N_\infty(V))$ consists of a unique zero function.
By Theorem~\ref{t6_3} this means that the moment problem is determinate.
$\Box$

Observe that the description of generalized resolvents for commuting symmetric operators (Theorem~\ref{t4_2})
can be used to study other moment problems, including truncated moment problems. The latter moment problems
lead to non-densely defined symmetric operators. Such moment and interpolation problems will be
studied elsewhere.

}

\noindent
Address:

V. N. Karazin Kharkiv National University \newline\indent
School of Mathematics and Computer Sciences \newline\indent
Department of Higher Mathematics and Informatics \newline\indent
Svobody Square 4, 61022, Kharkiv, Ukraine

Sergey.M.Zagorodnyuk@gmail.com; zagorodnyuk@karazin.ua

\end{document}